\documentclass{article}

\usepackage{hyperref}

\usepackage{amsfonts}
\usepackage{amssymb}
\usepackage{textcomp}
\usepackage{latexsym,amsmath}
\usepackage{graphics}

\renewcommand{\phi}{\varphi}
\newcommand{\be}{\begin{equation}}
\newcommand{\ee}{\end{equation}}
\newcommand{\ba}{\begin{eqnarray}}
\newcommand{\ea}{\end{eqnarray}}
\newcommand{\ban}{\begin{eqnarray*}}
\newcommand{\ean}{\end{eqnarray*}}

\newcommand{\nul}{{\bf0}}
\newcommand{\rd}{{\mathbb R}^d}
\newcommand{\zd}{{\mathbb Z}^d}
\newcommand{\td}{{\mathbb T}^d}
\renewcommand{\r}{{\mathbb R}}
\newcommand{\z} {{\mathbb Z}}
\newcommand{\cn} {{\mathbb C}}
\newcommand{\n} {{\mathbb N}}

\newcommand{\ddd}{,\dots,}
\renewcommand{\lll}{\left(}
\newcommand{\rrr}{\right)}
\newcommand{\h}{\widehat}
\newcommand{\w}{\widetilde}
\newcommand{\too}{\mathop{\longrightarrow}}

\title{Multivariate exact and falsified\\ sampling approximation
\thanks{This research was supported by Grants from RFBR 	(\#12-01-00216-a) and St.Petersburg State University (\#9.38.62.2012).}}

\author{
A. Krivoshein and M. Skopina
}
\date{Department of Applied Mathematics and Control Processes, 
 \\ St. Petersburg State University \\
 KrivosheinAV@gmail.com, skopina@ms1167.spb.edu}

\begin{document}
\maketitle

\begin{abstract}
Approximation properties of the expansions 
$\sum_{k\in\zd}c_k\phi(M^jx+k)$, where $M$ is a matrix dilation, 
$c_k$ is either the sampled value of a signal $f$ at $M^{-j}k$ or 
the integral average of $f$ near $M^{-j}k$ (falsified sampled value), are studied. 
Error estimations in $L_p$-norm, $2\le p\le\infty$, are given in terms of the Fourier transform of $f$. The approximation  order depends on how smooth is $f$, on the order of Strang-Fix condition for $\phi$ and on $M$. Some special properties of $\phi$ are
 required. To estimate the approximation order of falsified
sampling expansions we compare them with a differential expansions
$\sum_{k\in\,\zd}  Lf(M^{-j}\cdot)(-k)\phi(M^jx+k)$, where $L$ is an appropriate
differential operator. Some concrete functions $\phi$ applicable for implementations are constructed. In particular,  compactly supported splines and band-limited 
functions can be taken as $\phi$. Some of these functions provide expansions 
interpolating a signal at the points $M^{-j}k$.
\end{abstract}




\newtheorem{theo}{Theorem}
\newtheorem{lem}[theo]{Lemma}
\newtheorem {prop} [theo] {Proposition}
\newtheorem {coro} [theo] {Corollary}
\newtheorem {defi} [theo] {Definition}
\newtheorem {rem} [theo] {Remark}
\newtheorem {ex} [theo] {Example}

\newcommand{\tocsecindent}{\hspace{0mm}}

\section{Introduction}

The well-known sampling theorem (Kotel'nikov's or Shannon's formula)
states that
\be
f(x)=\sum_{k\in\z}  f(2^{-j}k)\,\frac{\sin\pi(2^jx-k)}{\pi (2^jx-k)}
\label{0}
\ee
for  band-limited to $[-2^{j-1},2^{j-1}]$ signals (functions) $f$.
 This formula is very useful for engineers.
It was just Kotel'nikov~\cite{02} and Shannon~\cite{01} who started to
apply this formula for signal processing, respectively in 1933 and 1949.
Up to now, an overwhelming diversity of digital signal processing applications
and devices are based on it and  more than successfully use it.
Without sampling theorem it would be impossible to make use of internet,
make photos and videos.
 However, mathematicians knew this formula much earlier, actually, it can be found
in the papers by Ogura~\cite{03} (1920), Whittaker~\cite{W} (1915),
Borel~\cite{Borel} (1897), and even Cauchy~\cite{04} (1841).

 Nowadays (\ref0) is also an important and
interesting formula for mathematicians.  Butzer with
co-authors recently  published several papers~\cite{Butz1}, \cite{Butz2},
\cite{Butz3},  where they analyze
sampling theorem, its applications and development.
In particular,  the equivalence
of sampling theorem to some other classical formulas was established
for a classes of band-limited functions. Also in~\cite{Butz6}, \cite{Butz7}
they studied a generalization of sampling decomposition replacing the
 sinc-function ${\rm sinc}(x):=\frac{\sin\pi x}{\pi x}$
 by certain  linear combinations of B-splines.

Equality~(\ref0)  holds only for functions $f\in L_2(\r)$
whose Fourier transform is supported on $[-2^{j-1},2^{j-1}]$.
However the right hand side of~(\ref0) (the sampling expansion of $f$)
has meaning for every continuous $f$ with a good enough decay.
The problem of approximation of $f$ by its sampling expansions as $j\to+\infty$
was studied by many mathematicians. We  mention only some of such results.
 Brown~\cite{Brown} proved that  for every $x\in\r$
\be
\left|f(x)- \sum_{k\in\,\z} f(-2^{-j}k)\,{\rm sinc}(2^jx+k)\right|\le
C \int\limits_{|\xi|>2^{j-1}}
 |\h f(\xi)|d\xi,
\label{brown}
\ee
whenever the Fourier transform of  $f$ is summable on $\r$.
It is known that the pointwise approximation by sampling
expansions does not hold for arbitrary continuous $f$, even compactly supported.
Moreover, Trynin~\cite{Tr4} proved that there exists a continuous function vanishing
outside of $(0,\pi)$ such that its  deviation from the sampling expansion diverges
at every point $x\in(0,\pi)$.  Approximation by sampling expansions  in $L_p$-norm
was  actively studied.
Bardaro,  Butzer,  Higgins, Stens and Vinti~\cite{Butz4}, \cite{Butz5},
 proved that
$$
\Delta_p:=\Big\|f- \sum_{k\in\,\z} f(-2^{-j}k)\,{\rm sinc}(2^j\cdot+k)\Big\|_p
\too\limits_{j\to+\infty}0\quad 1\le p<\infty,
$$
for $f\in C\cap \Lambda_p$, where $\Lambda_p$ consists of $f$
such that
$$
\sum_k|f(x_k)|^p(x_{k}-x_{k-1})<\infty
$$
for some class of admissible
partition $\{x_k\}_k$ of $\r$. Also they proved that the Sobolev spaces
${W_p^n}$, $n\in \n$, are subspaces of $\Lambda_p$, and that for every $f\in {W_p^n}$
\be
\Delta_p\le\frac{C\omega(f^{(n)}, 2^{-j})_p}{2^{-jn}}
\label{02}
\ee
where $\omega(\ )_p$ is the modulus of continuity in $L_p$.
The author of~\cite{Sk1}, \cite{Sk2} investigated approximation by sampling expansions
$$
\sum_{k\in\,\z} f(-2^{-j}k)\phi(2^jx+k)
$$
 for a wide class of band-limited functions $\phi$. For $p\ge2$,
 the error analysis was given in terms of the Fourier transform of $f$.
 In particular, the approximation order was found for functions $f$
 in Sobolev spaces $W_1^1$ with $f'\in Lip\ \alpha$, $\alpha>0$.
 In the case  $\alpha < 1-1/p$, the order of approximation is less
 than $2^{-j}$, which cannot be obtained from~(\ref{02}).
Similar results were proved for the generalized
sampling expansions (differential expansions)
\be
\sum_{k\in\,\z}  Lf(2^{-j}\cdot)(-k)\phi(2^jx+k),
\label{03}
\ee
where $Lf:=\sum_{l=0}^m{\alpha_l}f^{(l)}$. Also an analog of
Brown's estimate~(\ref{brown}) was prove for such expansions in~\cite{Sk2}.

Note that  differential  expansions~(\ref{03})
 can be useful for engineering applications. Indeed, engineers do not deal with
functions, they only have some discrete information about the function.
If values of  function at  equidistributed nodes are known, then
sampling expansion is very good for  recovering the function.
However, sometimes the values are known approximately.
 Assume that some
device gives the average value of a function $f$ on the interval
$[2^{-j}k, 2^{-j}(k+h)]$ instead of $f(2^{-j}k)$, i.e., one knows the values
\ban
\frac1{2^{-j}h}\int\limits_{2^{-j}k}^{2^{-j}(k+h)}f(t)\,dt\approx
\frac1{2^{-j}h}\int\limits_0^{2^{-j}h}\sum_{l=0}^m\frac1{l!}f^{(l)}(2^{-j}k)t^l\,dt=
\\
\sum_{l=0}^m\frac1{(l+1)!}h^l\frac{d^lf(2^{-j}\cdot)}{dx^l}(k).
\ean
But the latter sum is nothing as $Lf(2^{-j}\cdot)(k)$, where
$\alpha_l=\frac1{(l+1)!}h^l$.

Falsified sampling expansions (with the integral averages instead
of the exact sampled values), were studied by Butzer and Lei in~\cite{Butz8}.
They compared such generalized expansions with the corresponding usual
sampling expansions. Some error estimates  in the $L_\infty$-norm were presented.

The aim of the present paper is to extend the results of~\cite{Sk1}, \cite{Sk2}
to the multivariate case. However the class of functions $\phi$
we are considering
in the present paper is different.
We assume that $\phi$ is in  ${\cal L}_p$
$\h\phi$ has some appropriate decay.
(the space ${\cal L}_p$ is a subspace of $L_p$ introduced by Jia in~\cite{JiaMicPrewav}, see
Section~\ref{sa}).
The function $\rm sinc$ does not belong to this class, but the class includes
some compactly supported splines, which have an important advantage
in applied aspect because the corresponding sampling expansions are finite.
We study convergence and approximation order of  differential
sampling expansions in $L_p$-norm, $2\le p\le\infty$.
Also we analyze  error estimation for the deviation of 
falsified sampling expansions (with the integral average sampled values) 
from the corresponding differential  expansion.
It appeared that choosing compactly supported $\phi$ properly,
it is possible to provide almost the same order of the error as the
order of differential sampling  approximation. Our error analysis 
of falsified sampling expansions is new even for  the 
one-dimensional case.

The paper is organized as follows. Basic notations are given in
	Section~\ref{notation}. Some auxiliary lemmas are stated in
	Section~\ref{sa}. Section~\ref{scaleAppr} is devoted to
	the scaling approximation.
	Differential and falsified sampling approximation in $L_p$-norm 
		are discussed in Sections~\ref{GenSamExp} and~\ref{FalSamExp} respectively.
	In Section~\ref{example} some examples are given.

\section{Notations and basic facts}
\label{notation}

$\n$ is the set of positive integers,
    $\r$  is the set of real numbers,
    $\cn$ is the set of complex numbers.
    $\rd$ is the
    $d$-dimensional Euclidean space,  $x = (x_1\ddd x_d)$ and $y =
    (y_1\ddd y_d)$ are its elements (vectors),
    $(x)_j=x_j$ for $j~=~1,\dots,d,$
    $(x, y)~=~x_1y_1+~\dots~+x_dy_d$,
    $|x| = \sqrt {(x, x)}$, ${\bf0}=(0\ddd 0)\in\rd$;
      $B_r=\{x\in\rd:\ |x|\le r\}$, $\td=[-\frac 12,\frac 12]^d$;
    $\zd$ is the integer lattice
    in $\rd$, $\z_+^d:=\{x\in\zd:~x\geq~{\bf0}\}.$
    If $\alpha,\beta\in\zd_+$, $a,b\in\rd$, we set
    $[\alpha]=\sum\limits_{j=1}^d \alpha_j$,
    $\alpha!=\prod\limits_{j=1}^d\alpha_j!$,
    $\binom{\beta}{\alpha}=\frac{\alpha!}{\beta!(\alpha-\beta)!}$,
    $a^b=\prod\limits_{j=1}^d a_j^{b_j}$,
    $D^{\alpha}f=\frac{\partial^{[\alpha]} f}{\partial x^{\alpha}}=\frac{\partial^{[\alpha]} f}{\partial^{\alpha_1}x_1\dots
    \partial^{\alpha_d}x_d}$, 
    $\delta_{ab}$~is the Kronecker delta; if
     $n\in\n$, then
    $\Delta_n~:=~\{\alpha:~\alpha\in\zd_+,\,\,[\alpha]<n\}.$

An integer $d\times d$ matrix $M$ whose
eigenvalues are bigger than 1 in module is called a  dilation matrix.
Throughout the paper we
consider that such a matrix $M$ is fixed and  $m=|detM|$,
$M^*$ denote the conjugate matrix to $M$.
Since the spectrum of the operator $M^{-1}$ is
located in  $B_r$,
where $r=r(M^{-1}):=\lim\limits_{j\to+\infty}\|M^{-j}\|^{1/j}$
is the spectral radius of $M^{-1}$, and there exists at least
one point of the spectrum on the boundary of $B_r$, we have
	\be
	\|M^{-j}\|\le {C_{M,\theta}}\,\theta^{-j},\quad j\ge0,
	\label{00}
	\ee
for every  positive number $\theta$  which is smaller in module
than any eigenvalue of  $M$.
In particular, we can take $\theta > 1$, then
	\be
	\lim_{j\to+\infty}\|M^{-j}\|=0.
	\label{000}
	\ee
A  matrix $M$ is called isotropic
	if it is similar to a diagonal matrix
such that numbers $\lambda_1,\dots,\lambda_d$ are placed on the main diagonal
	and $|\lambda_1|=\cdots=|\lambda_d|$.
	Thus, $\lambda_1,\dots,\lambda_d$ are eigenvalues of $M$
	and the spectral radius of $M$ is equal to $|\lambda|,$
	where $\lambda$ is one of the eigenvalues of $M.$
	Note that if matrix $M$ is isotropic then
	$M^*$ is isotropic and $M^j$ is isotropic for all $j\in\z.$	
	It is well known that for an isotropic matrices $M$ and for  any $j\in\z$ 
	we have
	\be
	C^M_1 |\lambda|^j \le \|M^j\| \le C^M_2 |\lambda|^j,
	\label{10}
\ee
   where $\lambda$ is one of the eigenvalues of $M.$

If $\phi$ is a function defined on $\rd$, we set
$$
\phi_{jk}(x):=m^{j/2}\phi(M^jx+k),\quad j\in\z, k\in\rd.
$$

$L_p$ denotes $L_p(\rd)$, $1\le p\le\infty$. We use $W_p^n$, $1\le p\le\infty$, $n\in\n$, to denote
 the Sobolev space on~$\rd$, i.e. the set of
functions whose derivatives up to order $n$ are in $L_p(\rd)$,
with usual Sobolev norm.

If $f, g$ are functions defined on $\rd$ and $f\overline g\in L_1$,
then  $\langle  f, g\rangle:=\int\limits_{\rd}f\overline g$.

 For any function $f$, we set $f^{-}(\xi):=f(-\xi).$

If $F$ is a $1$-periodic (with respect to each variable) function and
$F\in L_1(\td)$,  then  $\widehat
F(k)=\int_{\td} F(x)e^{-2\pi i
(k,\,x)}\,dx$ is its $k$-th Fourier coefficient.
If $f\in L_1$,  then its Fourier transform is $\widehat
f(\xi)=\int_{\rd} f(x)e^{-2\pi i
(x,\xi)}\,dx$. 

Denote by $S$ the Schwartz class of functions defined on $\rd$.
    The dual space of $S$ is $S'$, i.e. $S'$ is
    the space of tempered distributions.
    The basic facts from distribution theory
    can be found, e.g., in~\cite{Vladimirov-1}.
    Suppose $f\in S$, $\phi \in S'$, then     
    $\langle \phi, f\rangle:= \overline{\langle f, \phi\rangle}:=\phi(f)$.
    If  $\phi\in S',$  then $\h \phi$ denotes its  Fourier transform
    defined by $\langle \h f, \h \phi\rangle=\langle f, \phi\rangle$,
    $f\in S$.
    If  $\phi\in S'$, $j\in\z, k\in\zd$, then we define $\phi_{jk}$ by
        $
        \langle f, \phi_{jk}\rangle=
        \langle f_{-j,-M^{-j}k},\phi\rangle, \forall f\in S.
        $

\section{Auxiliary results}
\label{sa}

Here and in what follows M denotes  a dilation matrix.
Given $\delta>0$, we introduce a special notation for
the following integrals if they make sense
$$
{\cal I}_{j,\gamma,q}^{In}(g)={\cal I}_{j,\gamma,q}^{In}(g, M, \delta)=
\int\limits_{|M^{*-j}\xi|<\delta}
	|\xi|^{q\gamma}|g(\xi)|^q d\xi, \quad
	{\cal I}_{j,\gamma,q}^{Out}(g)={\cal I}_{j,\gamma,q}^{Out}(g,M,\delta)=
\int\limits_{|M^{*-j}\xi|\ge\delta}
	|\xi|^{q\gamma}| g(\xi)|^q d\xi.
$$

\begin{lem}
\label{lem1}
Let $1\le q <\infty$, $1/p+1/q=1$, $j\in\z_+$,
	$\w\phi$ be a tempered distribution
	whose Fourier transform $\h{\w\phi}$ is a function on $\rd$
	such that $|\h{\w\phi}(\xi)|\le C_{\w\phi} |\xi|^{N}$
	 for almost all $\xi\notin\td$, $N=N({\w\phi})\ge 0,$ and
	 $|\h{\w\phi}(\xi)|\le C'_{{\w\phi}}$
	 for almost all $\xi\in\td$.
	 Suppose $g\in L_q$, $g(\xi)=O(|\xi|^{-N-d-\varepsilon})$
as $|\xi|\to\infty$, where $\varepsilon>0$;  $\gamma\in
(N+\frac dp, N+\frac dp+\epsilon)$
 for $q\ne1$, $\gamma=N$ for $q=1$, and set
$$
	G_j(\xi)=G_j({\w\phi},g,\xi):=\sum\limits_{l\in\,\zd}
	g(M^{*j}(\xi+l))\overline{\h{\w\phi}(\xi+l)}.
$$
 	Then $G_j$ is a $1$-periodic function in
  	$L_q(\td)$, $\langle g,\h{{\w\phi}_{jk}}\rangle=m^{j/2}\h G_j(k)$,
 and for every $\delta\in(0,\frac 12)$
 	 \be
	\left\|G_j-g(M^{*j}\cdot)\h{\w\phi}\right\|^q_{L_q(\td)}\le
    m^{- j} (C_{\gamma,\,{\w\phi}})^q \|M^{*-j}\|^{\gamma q} \,
   {\cal I}_{j,\gamma,q}^{Out}(g).
	\label{fLem1GjLq}
 	\ee
 	
\end{lem}

{\bf Proof.} First we will  prove~(\ref{fLem1GjLq}). For  $q>1$,
using  H\"older's inequality, we have
$$
	\left\|G_j-g(M^{*j}\cdot)\h{\w\phi}\right\|^q_{L_q(\td)}\le
    C^q_{\w\phi}\int\limits_{\td}
    \left|\sum\limits_{l\in\,\zd,\, l\neq \nul}
    \frac{|\xi+l|^\gamma| g(M^{*j}(\xi+l))|}{|\xi+l|^{\gamma-N}}
    \right|^q d\xi\le
$$
$$
    C^q_{\w\phi}
    \sup_{\xi\in\td}
    \left( \sum\limits_{l\in\,\zd,\, l\neq \nul}
    \frac{1}{|\xi+l|^{(\gamma-N)p}}\right)^\frac qp
 \int\limits_{\td}
        \sum\limits_{l\in\zd,\, l\neq \nul}
        	|(\xi+l)|^{\gamma q} |g(M^{*j}(\xi+l))|^q  d\xi =
           	$$
 $$
    (C_{\gamma,\,{\w\phi}})^q  \int\limits_{\rd\setminus\td}
                	|\xi|^{\gamma q} |g(M^{*j}\xi)|^q  d\xi \le
  (C_{\gamma,\,{\w\phi}})^q  \int\limits_{|\xi|\ge\delta}
                	|\xi|^{\gamma q} |g(M^{*j}\xi)|^q  d\xi \le
           	$$
           $$
            m^{-j} (C_{\gamma,\,{\w\phi}})^q \|M^{*-j}\|^{\gamma q} \,
    {\cal I}_{j,\gamma,q}^{Out}(g).
           $$
If $q=1$, then
$$
	\left\|G_j-g(M^{*j}\cdot)\h{\w\phi}\right\|_{L_1(\td)}\le
    C_{\w\phi}\int\limits_{\td}
    \sum\limits_{l\in\,\zd,\, l\neq \nul}
    {|\xi+l|^N| g(M^{*j}(\xi+l))|}      d\xi=
$$
 $$
        C_{\w\phi} \int\limits_{\rd\setminus\td}
                	|\xi|^{\gamma } |g(M^{*j}\xi)|  d\xi \le
      C_{\w\phi}  \int\limits_{|\xi|\ge\delta}
                	|\xi|^{\gamma } |g(M^{*j}\xi)|  d\xi \le
                      m^{-j}     C_{\w\phi}\|M^{*-j}\|^{\gamma } \,
    {\cal I}_{j,\gamma,1}^{Out}(g).
           $$
Combining~(\ref{fLem1GjLq}) with
$$
\int\limits_{\td}|
	g(M^{*j}\xi)\overline{\h{\w\phi}(\xi)}|^q d\xi\le
	(C'_{\w\phi})^q \int\limits_{\td}
	|g(M^{*j}\xi)|^q d\xi\le (C'_{\w\phi})^q m^{-j} \|g\|^q_q,
$$
using Minkowski's inequality and taking into account that
the integral ${\cal I}_{j,\gamma,q}^{Out}(g)$ is convergent,
we conclude that $G_j\in L_q(\td)$.

It remains to show that $\langle g,\h{{\w\phi}_{jk}}\rangle=m^{j/2}\h G_j(k)$,
$k\in\zd$,  but this is true because
    \ban
    \langle g,\h{{\w\phi}_{jk}}\rangle
    =m^{-j/2}\int\limits_{\rd} g(\xi)\overline{\h{\w\phi}(M^{*-j}\xi)}
    e^{-2\pi i\,(k,\,M^{*-j}\xi)}d\xi=
    m^{j/2}\int\limits_{\td}
    G_j(\xi) e^{-2\pi i\,(k,\,\xi)}d\xi=m^{j/2}\h G_j(k),
   \quad
    \ean
    where the last but one equality is justified by
    Lebesgue's dominated convergence theorem with taking into account
    the summability of the dominating function $\sum_{l\in\zd}
	|g(M^{*j}(\xi+l))\overline{\h{\w\phi}(\xi+l)}|$ on~$\td$,
	which follows from the proof 	of~(\ref{fLem1GjLq})
for $q=1.$
    $\Diamond$

Let $1\le p \le \infty$. Denote by ${\cal L}_p$ the set
	$$
	{\cal L}_p:=
	\left\{
	f\in L_p(\rd): \|f\|_{{\cal L}_p}:=
	\left\|\sum_{k\in\zd} \left|f(\cdot+k)\right|\right\|_{L_p(\td)}<\infty
	\right\}.
	$$
	With the norm $\|\cdot\|_{{\cal L}_p}$, ${\cal L}_p$ is a Banach space.
	The simple properties are:
	${\cal L}_1=L_1,$
	$\|f\|_p\le \|f\|_{{\cal L}_p}$,
	$\|f\|_{{\cal L}_q}\le \|f\|_{{\cal L}_p}$
	for $1\le q \le p \le\infty.$ Therefore, ${\cal L}_p\subset L_p$
	and ${\cal L}_p\subset {\cal L}_q$ for $1\le q \le p \le\infty.$
	If $f\in L_p$ and compactly supported then $f\in {\cal L}_p$ for $p\ge1.$
	If $f$ decays fast enough, i.e. there exist constants $C>0$
	and $\varepsilon>0$ such that
	$|f(x)|\le C( 1+|x|)^{-d-\varepsilon}$ $ \forall x\in\rd,$
	then $f\in {\cal L}_\infty$.
	Also we need the following statement
	(see~\cite[Theorem 2.1]{JiaMicPrewav})
	\begin{prop}
\label{propLp}
	Let $1\le p \le \infty$. If $\phi \in {\cal L}_p(\rd)$ and
	a sequence $a$ is in $\ell_p(\zd)$, then
	$$
\left\|\sum_{k\in\zd} a_k \phi_{0k}\right\|_p\le
	\|\phi\|_{{\cal L}_p} \|a\|_{\ell_p}.
$$
\end{prop}

	\begin{lem}
\label{lemQjLp}
	Let $g$ and $\w\phi$  be as in Lemma~\ref{lem1},
	Suppose $2\le p \le\infty$, $1/p+1/q=1$,
	 $\phi \in {\cal L}_p$.
	Then the series
	$\sum_{k\in\zd} \langle g,\h{\w\phi}_{jk}\rangle \phi_{jk}$
	converges unconditionally in $L_p(\rd)$.

\end{lem}

{\bf Proof.}
Suppose $2\le p <\infty$.
Because of Lemma~\ref{lem1} and the  Hausdorf-Young inequality, we
have
	$$
	\left(\sum_{k\in\zd} |\langle g,\h{{\w\phi}_{jk}}\rangle|^p\right)^\frac 1p =
	m^{\frac j2}
	\left(\sum_{k\in\zd} |\h G_j(k)|^p\right)^\frac 1p \le
	m^{\frac j2}  \|G_j\|_q <\infty,
	$$
where $G_j$ is a function from Lemma~\ref{lem1}.
By Proposition~\ref{propLp},  we can state that for every
finite	subset $\Omega$ of $\zd$
	$$\left\|\sum_{k\in\Omega}
	 \langle g,\h{{\w\phi}_{jk}}\rangle \phi_{jk}\right\|_p=
m^{\frac j2-\frac jp}\left\|\sum_{k\in\Omega}
	 \langle g,\h{{\w\phi}_{jk}}\rangle \phi_{0k}\right\|_p
	\le
	m^{\frac j2-\frac jp}\|\phi\|_{{\cal L}_p}
	\left(\sum_{k\in\Omega}
	|\langle g,\h{{\w\phi}_{jk}}\rangle|^p\right)^\frac 1p.
$$
	The series $\sum_{k\in\zd} |\langle g,\h{{\w\phi}_{jk}}\rangle|^p$ is
	convergent, which yields that
	$\sum_{k\in\zd} \langle g,\h{{\w\phi}_{jk}}\rangle \phi_{jk}$
	converges unconditionally.

Similarly, if  $p=\infty$, then
 for every finite	subset $\Omega$ of $\zd$
$$
\left\|\sum_{k\in\Omega}
	 \langle g,\h{{\w\phi}_{jk}}\rangle \phi_{jk}\right\|_\infty=
m^{\frac j2}\left\|\sum_{k\in\Omega}
	 \langle g,\h{{\w\phi}_{jk}}\rangle \phi_{0k}\right\|_\infty
	\le
	m^{\frac j2}
	\|\phi\|_{{\cal L}_\infty}	
	\sup\limits_{k\in\Omega} |\langle g,\h{{\w\phi}_{jk}}\rangle|=
m^{ j}
	\|\phi\|_{{\cal L}_\infty}	
	\sup\limits_{k\in\Omega} |\h G_j(k)|,
$$
and the unconditional convergence follows  from the Riemann-Lebesgue theorem.
$\Diamond$

\section{Scaling Approximation }
\label{scaleAppr}

Scaling operator
$\sum_{k\in\zd} \langle f, {\w\phi}_{jk}\rangle \phi_{jk}$
is a good tool of approximation for many  appropriate pairs of functions
$\phi, \w\phi$. We are interested in such operators, where $\w\phi$ is a
tempered distribution, e.g., the delta-function or a linear
combination of its derivatives. In this case the inner product
$\langle f, {\w\phi}_{jk}\rangle $ has meaning only for functions $f$ in $S'$.
To extend the class of functions $f$ one can  replace
$\langle f, {\w\phi}_{jk}\rangle $  by
$\langle \h f, \h{\w\phi}_{jk}\rangle$ or by $\langle g, \h{\w\phi}_{jk}\rangle$,
where $f=\h g^-$, and  in this case we set
$$
Q_j(\phi,\w\phi, f)=
 	\sum_{k\in\zd} \langle g, \h{\w\phi}_{jk}\rangle \phi_{jk},\quad j\in\z_+.
$$
In this section we study approximation properties of such operators
for a large class of distributions~$\w\phi$.

A function $f\in L_1(\rd)$ is said to satisfy the
{\em Strang-Fix condition of order $n$}
if  $D^{\beta} \h f (k)=0,$ for all $k\in \zd\setminus\{\nul\}$
and $\beta\in\Delta_n.$

	\begin{theo}
\label{theoQj}
	Let $2\le p \le \infty$, $1/p+1/q=1$.
	Suppose
\begin{itemize}
  \setlength{\itemsep}{0cm}%
  \setlength{\parskip}{0cm}%

 	\item[$\star$] $\phi \in {\cal L}_p $ and
	there exists $B_{\phi}>0$ such that
\be
\sum\limits_{k\in\zd}  |\h\phi(\xi+k)|^q<B_{\phi}\quad\forall \xi\in\r,
\label{2}
     \ee
      \item[$\star$]  $\w\phi$ be a tempered distribution
	 whose Fourier transform $\h{\w\phi}$ is a function on $\rd$
	such that $|\h{\w\phi}(\xi)|\le C_{\w\phi} |\xi|^{N}$
	 for almost any $\xi\notin\td$, $N=N(\w\phi)\ge 0,$ and
	 $|\h{\w\phi}(\xi)|\le C'_{{\w\phi}}$
	 for almost all $\xi\in\td$;
	
      \item[$\star$]  $f=\h g^{-}$, where
       $g\in L_q$, $g(\xi)=O(|\xi|^{-N-d-\varepsilon})$
as $|\xi|\to\infty$,  $\varepsilon>0$.

\end{itemize}

\noindent
 		Then
\begin{itemize}
	\item
	if  $\h{\w\phi}$ is continuous at zero
	and $\h\phi(0)\h{\w\phi}(0) = 1$, then the Strang-Fix condition of
	order $1$ for $\phi$ is necessary and
	sufficient for the convergence of
	$Q_j(\phi,\w\phi, f )$ to $f$ in the $L_p$-norm as $j\to +\infty$;

	\item
	if there exist $n\in\n$ and $\delta\in(0, 1/2)$ such that
$\h\phi\h{\w\phi}$ is  boundedly differentiable up to order $n$ on
	$\{|\xi|<\delta\}$,
$\h\phi$ is boundedly differentiable up to order $n$  on
	 $\{|\xi+l|<\delta\}$ for all $l\in\zd\setminus\{\nul\}$;
the function
	$\sum\limits_{l\in\zd,\, l\neq\nul}|D^\beta \h \phi (\xi+l)|$
	is bounded on	$\{|\xi|<\delta\}$  for $[\beta]=n$;
$D^{\beta}(1-\h\phi\h{\w\phi})(0) = 0$
	for all $\beta\in\Delta_{n}$;
	the Strang-Fix condition of order $n$ holds 	for $\phi$;
	$\gamma\in(N+\frac dp, N+\frac dp+\epsilon)$ for $p\ne\infty$,
 and $\gamma=N$ for $p=\infty$, then
	\be
	\|f-Q_j(\phi,\w\phi, f)\|_p^q\le C_1
	 \|M^{*-j}\|^{\gamma q}  {\cal I}_{j,\gamma,q}^{Out}(g)+
	 C_2 \|M^{*-j}\|^{nq}   {\cal I}_{j,n,q}^{In}(g),
	 \label{fTheoQjMain}
	\ee
where $C_1$ and $C_2$ do not depend on $j$ and $f$.

	\end{itemize}
\end{theo}	
	
{\bf Proof.} First of all we note that $\h\phi\in L_q$ due to~(\ref{2}),
and $\h\phi$ is continuous because $\phi\in L_1$.
Let  $G_j(\xi)=G_j(\widetilde\phi, g,\xi)$ be defined as in Lemma~\ref{lem1}.
This function is 1-periodic and  in $L_q(\td)$.
The function $G_j(M^{*-j}\cdot)\h\phi(M^{*-j}\cdot)$ is in $L_q$ because,
due to~(\ref{2}),
$$
\int\limits_{\rd}|G_j(M^{*-j}\xi)\h\phi(M^{*-j}\xi)|^q\,d\xi=
m^j\int\limits_{\td}|G_j(\xi)|^q
\sum\limits_{k\in\zd}|\h\phi(\xi+k)|^q\,d\xi\le m^jB_\phi
\|G_j\|^q_{L_q(\td)}.
$$
Hence, its  Fourier transform is in $L_p$ due to the Hausdorf-Young inequality.
On the other hand,  $Q_j(f)\in L_p$ by Lemma~\ref{lemQjLp},
and, due to Carleson's theorem (with the convergence of  Fourier series
over the cubic partial sums),
$$
G_j(M^{*-j}\xi)\h\phi(M^{*-j}\xi)=\sum\limits_{k\in\zd}\h G_j(k)e^{2\pi i(k,M^{*-j}\xi)}
\h\phi(M^{*-j}\xi)= \sum\limits_{k\in\zd}\langle g,\h{{\w\phi}_{jk}}
\rangle\widehat{\phi_{jk}}(\xi)\ \ a.e.
$$
The latter function coincides with $ \h{Q_j(f)}$  as a tempered distribution.
It follows from the du Bois-Reymond lemma that the  Fourier transform
of $G_j(M^{*-j}\cdot)\h\phi(M^{*-j}\cdot)$
coincides with $ {Q_j}(f)^-$ almost everywhere.	
Applying the Hausdorf-Young inequality we obtain
     \ba
      \|f-Q_j(f)\|_p= \|f^{-}-Q_j(f)^{-}\|_p\le
      \left\|g - G_j(M^{*-j}\cdot)\h\phi(M^{*-j}\cdot)\right\|_q\le
     \nonumber
\\
     \left\|g (1 - \h \phi(M^{*-j} \cdot)
     \overline{\h {\w\phi}(M^{*-j} \cdot)})\right\|_q +
     \left\| \h\phi(M^{*-j} \cdot)
     \sum\limits_{l\in\,\zd, l\neq 0}
     g( \cdot +M^{*j}l)
     \overline{\h{\widetilde\phi}(M^{*-j} \cdot +l)}\right\|_q.
     \label{theoQjf1}
     \ea
Let us fix dome $\delta\in (0,\frac 12)$ and consider the first summand in~(\ref{theoQjf1}).
	Changing the variable we obtain
	\ba
	\left\|g (1 - \h \phi(M^{*-j} \cdot)
     \overline{\h {\w\phi}(M^{*-j} \cdot)})\right\|^q_q=
     m^{j}
     \left\|g(M^{*j} \cdot) (1 - \h \phi
     \overline{\h {\w\phi}})\right\|^q_q=
     \nonumber
\\
     m^{ j}
     \int\limits_{|\xi|\ge\delta}
     |g(M^{*j}\xi) (1 - \h \phi(\xi)
     \overline{\h {\w\phi}(\xi)})|^q d\xi+
      m^{ j}  \int\limits_{|\xi|<\delta}
     |g(M^{*j}\xi) (1 - \h \phi(\xi)
     \overline{\h {\w\phi}(\xi)})|^q d\xi
    =:I_1+I_2.
     \label{theoQjFirstNorm}
	\ea
	
 Since there exists $C_{\phi,\w\phi}$ and $C_{\phi,\w\phi}'$
such that	$|1 - \h \phi(\xi)
     \overline{\h {\w\phi}(\xi)}|\le C_{\phi,\w\phi} |\xi|^\gamma$
     for almost all      $\xi\notin\td$ and
     $|1 - \h \phi(\xi)
     \overline{\h {\w\phi}(\xi)}|\le C'_{\phi,\w\phi}$
      for almost all      $\xi\in\td$, then we have   	
	\ba
	 I_1     \le
	(C_{\phi,\w\phi})^q m^{j}
	\int\limits_{\rd\setminus\td}
	|\xi|^{\gamma q}|g(M^{*j}\xi)|^q d\xi+
\nonumber
\\
(C'_{\phi,\w\phi})^q\delta^{-\gamma q}m^{j}
	\int\limits_{\td\setminus\{|\xi|<\delta\}}
	|\xi|^{\gamma q}|g(M^{*j}\xi)|^q d\xi\le
	 C_{I_1}
	 \|M^{*-j}\|^{\gamma q} \,
	 {\cal I}_{j,\gamma,q}^{Out}(g).
	 \label{lem1RdTd}	
	\ea
	    Thus, $I_1\to 0$ as $j\to +\infty$.  The second integral $I_2$ is
     $$
     I_2=
      m^{ j}  \int\limits_{|\xi|<\delta}
     |g(M^{*j}\xi) (1 - \h \phi(\xi)
     \overline{\h {\w\phi}(\xi)})|^q d\xi =
      \int\limits_{\rd}
     \chi_{M^{*j}\{|\xi|<\delta\}}(\xi)
     |g(\xi) (1 - \h \phi(M^{*-j}\xi)
     \overline{\h {\w\phi}(M^{*-j}\xi)})|^q d\xi.
     $$
By the additional  assumptions of the first statement,
     the integrand tends to zero as $j\to +\infty$ for each $\xi\in\rd.$
     Since the integrand has a summable majorant, by
     Lebesque's  dominated convergence theorem,
     we conclude that  $I_2 \to 0$ as $j\to +\infty$.
     Thus, by~(\ref{theoQjFirstNorm}),
     \be
     \left\|g (1 - \h \phi(M^{*-j}  \cdot)
     \overline{\h {\w\phi}(M^{*-j}  \cdot)})
     \right\|_q \to 0\,\, as \,\,
     j\to +\infty.
     \label{theoQjFirstNormZero}
     \ee
Next, the second summand in~(\ref{theoQjf1}) is
	\ba
\left\| \h\phi(M^{*-j}  \cdot)\sum\limits_{l\in\,\zd, l\neq 0}
     g(  \cdot +M^{*j}l)
     \overline{\h{\widetilde\phi}(M^{*-j} \cdot +l)}\right\|_q^q
	= m^{j}    \left\| \h\phi\sum\limits_{l\in\,\zd, l\neq 0}
     g(M^{*j}( \cdot+l))
     \overline{\h{\widetilde\phi}( \cdot +l)}\right\|_q^q  =
    \nonumber
    \\
     m^{j}\sum\limits_{k\in\zd}
      \int\limits_{\td}
     |\h\phi(\xi+k)|^q
     \left|\sum\limits_{l\in\,\zd, l\neq 0}
     g(M^{*j}(\xi+l+k))
     \overline{\h{\widetilde\phi}(\xi+l+k)}\right|^q d\xi=
      \nonumber
    \\
  m^{j} \sum\limits_{k\in\zd}
      \int\limits_{\td}
     |\h\phi(\xi+k)|^q
     \left|\sum\limits_{l\in\,\zd, l\neq k}      g(M^{*j}(\xi+l))
     \overline{\h{\widetilde\phi}(\xi+l)}\right|^q d\xi=
  \nonumber
    \\
     m^{j}    \int\limits_{\td}
     |\h\phi(\xi)|^q
     \left|\sum\limits_{l\in\,\zd, l\neq 0}
     g(M^{*j}(\xi+l))
     \overline{\h{\widetilde\phi}(\xi+l)}\right|^q d\xi+
    \nonumber
    \\
     m^{j}  \int\limits_{\td}
      \sum\limits_{k\in\zd, k\neq\nul}
     |\h\phi(\xi+k)|^q
     \left|\sum\limits_{l\in\,\zd,  l\neq k}
     g(M^{*j}(\xi+l))
     \overline{\h{\widetilde\phi}(\xi+l)}\right|^q d\xi=:J_1+J_2.
     \label{3}
     \ea
     To estimate   $J_{1}$  we use~(\ref{2}), and then, by~(\ref{fLem1GjLq}),
     	\be
	J_1\le B_\phi m^j \left\|G_j-g(M^{*j}\cdot)\h{\w\phi}\right
\|^q_{L_q(\td)}\le C_{J_1}  \|M^{*-j}\|^{\gamma q}
    \,
	 {\cal I}_{j,\gamma,q}^{Out}(g).
	\label{4}
\ee
		Thus, $J_1\to 0$ as $j\to +\infty$.
	Using Minkowski's inequality,
the second summand $J_2$ can be estimated as follows
	\ba
	J_2^{1/q}\le
	\left(m^{j}
      \int\limits_{\td}
      \sum\limits_{k\in\zd, k\neq\nul}
     |\h\phi(\xi+k)|^q     \left| \sum\limits_{l\in\,\zd,  l\neq 0}
     g(M^{*j}(\xi+l))
     \overline{\h{\widetilde\phi}(\xi+l)}\right|^q d\xi
     \right)^\frac 1q +
       \nonumber
    \\
         \left( m^{j}
      \int\limits_{\td}
      \sum\limits_{k\in\zd, k\neq\nul}
     |\h\phi(\xi+k)|^q
     |g(M^{*j}(\xi+k))
     \overline{\h{\widetilde\phi}(\xi+k)}|^q d\xi
      \right)^\frac 1q +
    \nonumber
    \\
        \left( m^{j}
      \int\limits_{\td}
      \sum\limits_{k\in\zd, k\neq\nul}
     |\h\phi(\xi+k)|^q
     |g(M^{*j}\xi)
     \overline{\h{\widetilde\phi}(\xi)}|^q d\xi
      \right)^\frac 1q
     =:(J_{21})^{1/q}+(J_{22})^{1/q}+(J_{23})^{1/q}.
	\label{5}
\ea
To estimate   $J_{21}$ again we use~(\ref{2}), and then, by~(\ref{fLem1GjLq}),
     	\be
	J_{21}\le B_\phi m^j \left\|G_j-g(M^{*j}\cdot)\h{\w\phi}\right
\|^q_{L_q(\td)}\le C_{J_{21}}  \|M^{*-j}\|^{\gamma q}
    \, {\cal I}_{j,\gamma,q}^{Out}(g).
		\label{6}
\ee
Since $\h\phi$ is bounded,  similarly to~(\ref{lem1RdTd}), we have
\be
	J_{22}=m^{j}
      \int\limits_{\rd\setminus\td}  |\h\phi(\xi)
      g(M^{*j}\xi)
     \overline{\h{\widetilde\phi}(\xi)}|^q d\xi
      \le C_{J_{22}}  \|M^{*-j}\|^{\gamma q}
    \, {\cal I}_{j,\gamma,q}^{Out}(g).
	\label{7}
\ee	
	Thus, $J_{21}\to 0$, $J_{22}\to 0$, as $j\to +\infty$.
	

The third summand $J_{23}$ can be represented as
	     	$$
	J_{23}=
	 \int\limits_{\rd}
	\chi_{M^{*j}\td}(\xi)
      \sum\limits_{k\in\zd, k\neq\nul}
     |\h\phi(M^{*-j}\xi+k)|^q
     |g(\xi)
     \overline{\h{\widetilde\phi}(M^{*-j}\xi)}|^q d\xi
	$$
By the additional  assumptions of the first statement,
	 the integrand tends to    $\sum\limits_{k\in\zd, k\neq\nul}
     |\h\phi(k)|^q
     |g(\xi)\overline{\h{\widetilde\phi}(\nul)}|^q$ as
     $j\to +\infty$ for each $\xi\in\rd.$ Since the integrand
     has a summable majorant, by
     Lebesque's  dominated convergence theorem,
     we conclude that
     $$J_{23} \to \sum\limits_{k\in\zd, k\neq\nul}
     |\h\phi(k)|^q
     |\overline{\h{\widetilde\phi}(\nul)}|^q
     \|g\|_q^q,\,\quad j\to +\infty.$$
Thus,
     $$
     \left\| \h\phi(M^{*-j}  \cdot)
     \sum\limits_{l\in\,\zd, l\neq 0}
     g(  \cdot+M^{*j}l)
     \overline{\h{\widetilde\phi}(M^{*-j} \cdot+l)}\right\|^q_q
 \to      \sum\limits_{k\in\zd, k\neq\nul}
     |\h\phi(k)|^q
     |\overline{\h{\widetilde\phi}(\nul)}|^q
     \|g\|_q^q,\,\,\, as \,\,j\to +\infty,
     $$
     which together with~(\ref{theoQjf1})
     and~(\ref{theoQjFirstNormZero}) yields that
     the  Strang-Fix condition for $\phi$ is necessary and sufficient
	for $\|f-Q_j(f)\|_p \to 0$ as $j\to +\infty$.

Assume now that all assumptions of the second statement are satisfied. Set
$$
 B_{\phi,\w\phi,n}:=
     \sup\limits_{|\xi|<\delta}
     \left(\sum\limits_{[\beta]=n} \frac
     {|D^\beta [\h \phi \overline{\h{\widetilde\phi}}](\xi)|}
     {\beta!}\right)^q, \quad
      B_{\phi,n}= \sup\limits_{|\xi|<\delta}
     \left(
       \sum\limits_{l\in\zd,\, l\neq\nul}
     \left|\sum\limits_{[\beta]=n}
     \frac {D^\beta \h \phi (\xi+l)}
     {\beta!}\right|^q
     \right).
$$
     If $|\xi|<\delta,$ $l\neq \nul,$ then, by the Taylor formula
     with the remainder  in Lagrange's form, we get
     $$
     \h \phi(\xi+l)=
     \sum\limits_{[\beta]=n} \frac {\xi^\beta}{\beta!}
     D^\beta \h\phi(l+ t\xi) , \quad
     1- \h \phi(\xi) \overline{\h{\widetilde\phi}(\xi)}=
     -\sum\limits_{[\beta]=n} \frac {\xi^\beta}{\beta!}
     D^\beta [\h \phi \overline{\h{\widetilde\phi}}](r\xi) ,
     $$
     for some $t, r \in (0,1)$, and hence
     $$
     \sum_{l\in\zd,\, l\neq\nul}
     |\h \phi(\xi+l)|^q\le
     B_{\phi,n} |\xi|^{nq},\quad
     | 1- \h \phi(\xi) \overline{\h{\widetilde\phi}(\xi)}|^q\le
     B_{\phi,\w\phi,n} |\xi|^{nq}.
     $$
 It follows that
     $$
       m^{ j}  \int\limits_{|\xi|<\delta}
     |g(M^{*j}\xi) (1 - \h \phi(\xi)
     \overline{\h {\w\phi}(\xi)})|^q d\xi\le
     B_{\phi,\w\phi,n} \|M^{*-j}\|^{nq} \,    {\cal I}_{j,n,q}^{In}(g),
     $$
     and
      $$
       m^{ j}  \int\limits_{|\xi|<\delta}\sum\limits_{k\in\zd, k\neq\nul}
     |\h\phi(\xi+k)|^q
     |g(M^{*j}\xi)
     \overline{\h{\widetilde\phi}(\xi)}|^q \,d\xi\le
     B_{\phi,n}(C'_{\w\phi})^q
     \|M^{*-j}\|^{nq} \,    {\cal I}_{j,n,q}^{In}(g).
     $$
     Combining these relations
     with~(\ref{lem1RdTd}) and~(\ref{2})  respectively, we obtain
     $$
     m^{ j}  \int\limits_{\rd}
     |g(M^{*j}\xi) (1 - \h \phi(\xi)
     \overline{\h {\w\phi}(\xi)})|^q d\xi \le
           B_{\phi,\w\phi,n}\|M^{*-j}\|^{nq}   \,
            {\cal I}_{j,n,q}^{In}(g)
           +C_{I_1} \|M^{*-j}\|^{\gamma q} \,
           {\cal I}_{j,\gamma,q}^{Out}(g)
     $$
and
\ba
m^{j}\int\limits_{\td}\sum\limits_{k\in\zd, k\neq\nul}
     |\h\phi(\xi+k)|^q      |g(M^{*j}\xi)
     \overline{\h{\widetilde\phi}(\xi)}|^q d\xi\le
     B_{\phi,n}(C'_{\w\phi})^q  \|M^{*-j}\|^{nq} \,
     {\cal I}_{j,n,q}^{In}(g)+
     \nonumber
     \\
    B_{\phi}  m^{j}\int\limits_{\td\setminus \{|\xi|<\delta\}}
        |g(M^{*j}\xi)  \overline{\h{\widetilde\phi}(\xi)}|^q d\xi\le
       (C'_{\w\phi})^q  \,( B_{\phi,n}  \|M^{*-j}\|^{nq}
     {\cal I}_{j,n,q}^{In}(g)+
     B_{\phi} \delta^{-\gamma q}  \|M^{*-j}\|^{\gamma q}
       \, {\cal I}_{j,\gamma,q}^{Out}(g) ).
      \nonumber
  \ea
 These estimations  together with~(\ref{theoQjf1})-(\ref{lem1RdTd}) and
 (\ref{3})-(\ref{7})
 complete the proof of~(\ref{fTheoQjMain}).  $\Diamond$

 Note that the above proof is based on the technique
employed by Jetter and Zhou~\cite{JZ}, ~\cite{JZ1},
and developed in~\cite{KS}.

Unfortunately it is complicated to estimate
 the  approximation order of $Q_j(\phi,\w\phi,f)$
 from~(\ref{fTheoQjMain}) in Theorem~\ref{theoQj} for the general case.
 But we now will give such estimations for some special cases.

\begin{theo}
Let $M$ be an isotropic matrix dilation and $\lambda$ be its eigenvalue, $j\in\z_+$.
Suppose all conditions of Theorem~\ref{theoQj} are fulfilled.
 Then
	 $$\|f-Q_j(\phi,\w\phi,f)\|_p\le
	 \begin{cases}
	 C |\lambda|^{-j(N+\frac dp + \varepsilon)}  &\mbox{if }
	n> N+\frac dp + \varepsilon\\
	  C j^{1/q} |\lambda|^{-jn} &\mbox{if }
	 n= N+\frac dp + \varepsilon \\
	C|\lambda|^{-jn}
	 &\mbox{if }
	 n< N+\frac dp + \varepsilon
	\end{cases},$$
	where  $C$ does not depends on $j$.
	\label{theoQjOrder}
\end{theo}
	
{\bf Proof.} Throughout the proof we denote by $C$ and $C'$
different constants which do not depend on $j$.
Since $g(\xi) = O (|\xi|^{-N-d-\varepsilon})$,
there exists a big enough number $A\in\r$ such that
$|g(\xi)|\le C |\xi|^{-N-d-\varepsilon}$ for any $|\xi|>A.$
 By Theorem~\ref{theoQj} inequality~(\ref{fTheoQjMain}) is valid.
Let us consider the first term in~(\ref{fTheoQjMain}). Since
the set $\{|M^{*-j}\xi| \ge\delta\}$
	is a subset of
	$\{|\xi|\ge\delta/\|M^{*-j}\|\}$, we have
	$$
	\|M^{*-j}\|^{\gamma q}  {\cal I}_{j,\gamma,q}^{Out}(g)=
	 \|M^{*-j}\|^{\gamma q}
	\int\limits_{|M^{*-j}\xi|\ge\delta}
	|\xi|^{\gamma q}| g(\xi)|^q d\xi	\le
	C\|M^{*-j}\|^{\gamma q}
	\int\limits_{|\xi|\ge\delta/\|M^{*-j}\|}
	\frac {d\xi}
	{|\xi|^{(N+d+\varepsilon-\gamma)q}} ,
	$$
for all $j>j_0,$ where $j_0\in\z$ is such that 	$\frac {\delta} {\|M^{*-j_0}\|}>A.$
Using general polar coordinates with $\rho:=|\xi|$ and taking into account that
$(N+d+\varepsilon-\gamma)q>d$, we obtain
$$
\int\limits_{|\xi|\ge\delta/\|M^{*-j}\|}
	\frac {d\xi}
	{|\xi|^{(N+d+\varepsilon-\gamma)q}}\le
C\int\limits_{ {\delta}/{\|M^{*-j}\|}}^{+\infty}
	\frac {1}
	{\rho^{(N+d+\varepsilon-\gamma)q-d+1}} d\rho\le
C'\|M^{*-j}\|^{q(N+\frac dp + \varepsilon-\gamma)},
$$
and, by~(\ref{10}),
	\ba
	\|M^{*-j}\|^{\gamma q}  {\cal I}_{j,\gamma,q}^{Out}(g) \le	
C\|M^{*-j}\|^{q(N+\frac dp + \varepsilon)}\le
	C' |\lambda|^{-j  q(N+\frac dp + \varepsilon)}.
	\label{fTheoJ_estimate}
	\ea

Next, let us consider the second term in~(\ref{fTheoQjMain}).
  Since the set $\{|M^{*-j}\xi|<\delta\}$ is a subset of
  $\{|\xi|<\|M^{*j}\|\delta\}$, we get
  $$
  	\|M^{*-j}\|^{nq}   {\cal I}_{j,n,q}^{In}(g)=
	\|M^{*-j}\|^{n q}
	\int\limits_{|M^{*-j}\xi|<\delta}
	|\xi|^{n q}| g(\xi)|^q d\xi \le
  	C\|M^{*-j}\|^{n q}
	\int\limits_{|\xi|<\|M^{*j}\|\delta}
	|\xi|^{n q}| g(\xi)|^q d\xi \le
	$$
	$$
	C' \lll\|M^{*-j}\|^{nq}  +
	\|M^{*-j}\|^{n q}
	\int\limits_{A<|\xi|<\|M^{*j}\|\delta}
	|\xi|^{(n-N-d-\varepsilon) q}d\xi\rrr
  $$
    for all $j>j_1$, where $j_1\in\z$ is such that $A<\|M^{*j_1}\|\delta$.
Using general polar coordinates with $\rho:=|\xi|$ we obtain
	$$\|M^{*-j}\|^{nq}   {\cal I}_{j,n,q}^{In}(g) \le
	C \lll\|M^{*-j}\|^{nq}  +
	\|M^{*-j}\|^{n q}
	\int\limits_{A}^{\|M^{*j}\|\delta}
	\rho^{(n-N-d-\varepsilon) q+d-1}d \rho\rrr.
	$$
Let  $J:=\max\{j_0,j_1\}$, $j>J$. If $n< N+\frac dp + \varepsilon$, then
the integral
$$
\int\limits_{A}^{+\infty}
	\rho^{(n-N-d-\varepsilon) q+d-1}d \rho
$$
is convergent, and using~(\ref{10}), we have
	$
	\|M^{*-j}\|^{nq}   {\cal I}_{j,n,q}^{In}(g) \le
	C |\lambda|^{-jnq}.
	$
	Therefore,
	together with~(\ref{fTheoJ_estimate})  the final estimate in this case is
     $
	 \|f-Q_j(\phi,\w\phi,f)\|_p\le
	  C |\lambda|^{-jn }.
	 $

Let  $n=N+\frac dp + \varepsilon$. Again using~(\ref{10}), we have
	 $$\|M^{*-j}\|^{nq}   {\cal I}_{j,n,q}^{In}(g) \le
	 C\lll \|M^{*-j}\|^{nq} +  \|M^{*-j}\|^{nq}
	 \ln \|M^{*j}\|\rrr \le C' j |\lambda|^{-jnq}.
	 $$
	Together with~(\ref{fTheoJ_estimate})  the final estimate in this case is
	 $
	 \|f-Q_j(\phi,\w\phi,f)\|_p\le
	 C j^{1/q} |\lambda|^{-jn}.
	 $	
	
Similarly, if $n> N+\frac dp + \varepsilon$, then
	 		 $$\|M^{*-j}\|^{nq}   {\cal I}_{j,n,q}^{In}(g) \le
	 C\lll|\lambda|^{-jnq}  + |\lambda|^{-jq(N+\frac dp+\varepsilon)}\rrr\le
C' |\lambda|^{-jq(N+\frac dp+\varepsilon)},
	 $$	
Therefore, together with~(\ref{fTheoJ_estimate})  the final estimate in this case is
	 	$
	 \|f-Q_j(\phi,\w\phi,f)\|_p\le
	  C |\lambda|^{-j(N+\frac dp+\varepsilon)}.\ \Diamond $

Observing the proof of Theorem~\ref{theoQjOrder} for the case
$n<N+\frac dp + \varepsilon$, we see that the estimate does not depend
on $\|M^{*j}\|$. So,  in this case, we can repeat the proof for arbitrary
matrix dilation $M$ with using~(\ref{00})  instead of~(\ref{10}),
which leads to the following statement.

\begin{theo}
Suppose all conditions of Theorem~\ref{theoQj} are fulfilled,
and $n<N+\frac dp + \varepsilon$.
	Then
 \be
 \|f-Q_j(\phi,\w\phi, f)\|_p\le C \theta^{-jn}
 \label{9}
 \ee
 for every  positive number $\theta$   which is smaller in module
than any eigenvalue of  $M$ and  some $C$ which does not depend on $j$.
		\label{theoQjOrder1}
\end{theo}

The latter theorem  does not provide
 approximation order of $Q_j(\phi,\w\phi, f)$ better than $\theta^{-jn}$
 even for  very  smooth functions $f$.
This drawback can be fixed
under stronger restrictions on  $\phi$.

 	\begin{theo}
\label{theoQj1}
	Let $2\le p \le \infty$, $1/p+1/q=1$.
	Suppose
\begin{itemize}
  \setlength{\itemsep}{0cm}%
  \setlength{\parskip}{0cm}%

 	  \item[$\star$]  $\w\phi$ be a tempered distribution
	 whose Fourier transform $\h{\w\phi}$ is a function on $\rd$
	such that $|\h{\w\phi}(\xi)|\le C_{\w\phi} |\xi|^{N_{\w\phi}}$
	 for almost any $\xi\notin\td$, $N_{\w\phi}>0,$ and
	 $|\h{\w\phi}(\xi)|\le C'_{{\w\phi}}$
	 for almost all $\xi\in\td$;

\item[$\star$] $\phi \in {\cal L}_p $, the function
 $\sum_{k\in\zd}  |\h\phi(\xi+k)|^q$	is bounded, and
     there exists $\delta\in(0,1/2)$ such that
 $\overline{\h\phi}(\xi)\h{\w\phi}(\xi)=1$
 a.e. on $\{|\xi|<\delta\}$, $\h\phi(\xi)=0$ a.e.
 on $\{|l-\xi|<\delta\}$ for all
 $l\in\z\setminus\{0\}$;

      \item[$\star$]  $f=\h g^{-}$, where
       $g\in L_q$, $g(\xi)=O(|\xi|^{-N-d-\varepsilon})$
as $|\xi|\to\infty$,  $\varepsilon>0$.

\end{itemize}

\noindent
	Then
 \be
 \|f-Q_j(\phi,\w\phi, f)\|_p\le C \theta^{-j(N+\frac{d}{p}+\varepsilon)}
 \label{91}
 \ee
 for every  positive number $\theta$   which is smaller in module
than any eigenvalue of  $M$ and  some $C$ which does not depend on $j$.


\end{theo}

{\bf Proof.}
 Let
$\gamma\in(N+\frac dp, N+\frac dp+\epsilon)$ for $p\ne\infty$,
 and $\gamma=N$ for $p=\infty$.	
First we prove that
\be
	\|f-Q_j(\phi,\w\phi, f)\|_p^q\le C'
	 \|M^{*-j}\|^{\gamma q}  {\cal I}_{j,\gamma,q}^{Out}(g),
	 \label{8}
	\ee
where $C'$ does not depend on $j$ and $f$.
 For any compact set $K\subset \rd$,  function $g$ can be approximated
 in $L_q(K)$ by infinitely smooth functions  supported on $K$.
So, given $j$, one can find a function $\h F_j\in C^\infty(\rd)$ such that
${\rm supp}\,\h F_j\subset \{|M^{-j}\xi|<\delta\}$ and
\be
\int\limits_{|M^{-j}\xi|<\delta}|
g(\xi)-\h F_j(\xi)|^qd\xi\le
\left(\sup\limits_{|M^{-j}\xi|<\delta}|\xi|^{Nq}\right)^{-1}
 \|M^{*-j}\|^{(\gamma-N) q}  {\cal I}_{j,\gamma,q}^{Out}(g).
\label{35}
\ee
This yields
\ba
 {\cal I}_{j,N,q}^{In}(g-\h F)=
 \int\limits_{|M^{-j}\xi|<\delta}
|\xi|^{Nq}|(g(\xi)-\h F_j(\xi))|^qd\xi\le
\|M^{*-j}\|^{(\gamma-N) q}  {\cal I}_{j,\gamma,q}^{Out}(g).
\label{36}
\ea
Let $F_j$ be a function whose Fourier transform is $\h F_j$.
Evidently $F_j^-={\h{\h F}}_j$, and, due to  Carleson's theorem
and Lemma~\ref{lem1},  we have
$$
\sum\limits_{k\in\zd}\langle \h F,\h{{\w\phi}_{jk}}
\rangle\widehat{\phi_{jk}}(\xi)=
\sum\limits_{l\in\,\zd}
\h F_j(\xi+M^{j}l)\overline{\h{\w\phi}(M^{-i}\xi+l)}\h{\phi}(M^{-j}\xi)
$$
If $l\ne0$ and $F_j(\xi+M^{j}l)\ne0$, then $|M^{-j}\xi+l|<\delta$ and hence
$\h{\phi}(M^{-j}\xi)=0$. So,
$$
\sum\limits_{k\in\zd}\langle \h F_j,\h{{\w\phi}_{jk}}
\rangle\widehat{\phi_{jk}}(\xi)=\h F_j(\xi),
$$
which yields that $Q_j(\phi,\w\phi, F_j)=F_j$.
It follows that
$$
\|f- Q_j(\phi,\w\phi, f)\|_p=\|f-F_j- Q_j(\phi,\w\phi, f-F_j)\|_p.
$$
Since the assumptions of the second statement in Theorem~\ref{theoQj}
is satisfied for every $n$, in particular, for $n=N$,  we obtain
$$
\|f-F_j- Q_j(\phi,\w\phi, f-F_j)\|_p\le C_1
\|M^{*-j}\|^{N q}  {\cal I}_{j,N,q}^{In}(g-\h F)
+C_2 \|M^{*-j}\|^{\gamma q}  {\cal I}_{j,\gamma,q}^{Out}(g-\h F).
$$
To complete the proof of~(\ref{8}) it remains to use~(\ref{36})
and take into account that
$ {\cal I}_{j,\gamma,q}^{Out}(g-\h F)= {\cal I}_{j,\gamma,q}^{Out}(g)$.

Using the first inequality in~(\ref{fTheoJ_estimate}), we have
$$
\|M^{*-j}\|^{\gamma q}  {\cal I}_{j,\gamma,q}^{Out}(g) \le	
C_3\|M^{*-j}\|^{q(N+\frac dp + \varepsilon)},
$$
which yields~(\ref{9}) due to~(\ref{00}). 
 $\Diamond$

\section{Differential expansions}
\label{GenSamExp}	
	
The above  theorems  are proved for a wide
class of distributions $\w\phi$.
The $\delta$-function can be taken as $\w\phi$   with $N=0$.
If $f=\h g^-$ and $g$ as in Theorem~\ref{theoQj}, then $g\in L_1$  and
$$
f(-M^{-j}k)=\h g(M^{-j}k)=
\int\limits_{\rd} g(\xi)
	e^{-2\pi i (k,M^{*-j}\xi)}\,d\xi=m^{j/2} \langle g, \h\delta_{jk}\rangle,
\quad k\in\zd.
$$
Hence
$$	
Q_j(\phi,\w\phi, f)=m^{-j/2}\sum_{k\in\zd}f(-M^{-j}k)\phi_{jk},
$$	
i.e. $Q_j$ is a sampling expansion of $f$ in this case.	

Consider next a differential operator $L$ 
	defined by
	\be
Lf:=\sum_{\beta\in\Delta_{N+1}} a_\beta D^{\beta} f,
\quad a_{\beta}\in{\mathbb C}, a_{\nul}\ne0,
\label{11}
\ee
where	$N\in\z_+.$
	The action of operator $D^{\beta}$ is associated with the action of the corresponding derivative of the $\delta$-function.
	In  more detail, let
	$f=\h g^-$,
	$
	\int\limits_{\rd} (1+|\xi|)^{N+\alpha} |g(\xi)| d\xi < \infty,
	$
	$\alpha>0.$
	Therefore,  $f$ is  continuously differentiable on $\rd$
	up to the order $N.$
	Then for  $\beta\in\Delta_{N+1} $
	\ban
	D^{\beta}f (M^{-j}\cdot)(-k)=
	 D^{\beta}\h g (M^{-j}\cdot)(k)=
	(-1)^{[\beta]}m^j D^{\beta}\h {g(M^{*j}\cdot)} (k)=
\\
	(-1)^{[\beta]} m^j \int\limits_{\rd} g(M^{*j}\xi) (-2\pi i \xi)^{\beta}
	e^{-2\pi i (k,\xi)} d\xi=
	(-1)^{[\beta]} \int\limits_{\rd} g(\xi) \overline{
	(2\pi i M^{*-j}\xi)^{\beta} e^{2\pi i (M^{-j}k,\xi)}} d\xi=
\\
	(-1)^{[\beta]} m^{j/2}
	\langle
	g, \h {D^{\beta} \delta_{jk}}
	\rangle
\ean
	If now $\w\phi = \sum\limits_{\beta\in\Delta_{N+1}}
\overline{a_{\beta}} (-1)^{[\beta]} D^{\beta} \delta$ (we say that $\w\phi$
\textit{is associated with} $L$), then
	$$
m^{-j/2} L f(M^{-j}\cdot)(-k)=	
	\langle 	g, \h {\w\phi_{jk}} 	\rangle,\quad k\in\zd.
$$
Hence
$$
Q_j(\phi,\w\phi, f)=m^{-j/2}\sum_{k\in\zd}L f(M^{-j}\cdot)(-k)\phi_{jk}.
$$
	Thus,
	choosing  appropriate  function $\phi$ we can provide all
conditions of Theorem~\ref{theoQj}
which together with Theorems~\ref{theoQjOrder} and~\ref{theoQjOrder1}  gives
the following statement.
	
		\begin{theo}
\label{theoQjL}
	Let $2\le p \le \infty$, $1/p+1/q=1$, $N\in\z_+,$
$\gamma\in(N+\frac dp, N+\frac dp+\epsilon)$ for $p\ne\infty$,
 and $\gamma=N$ for $p=\infty$, a differential operator $L$ be defined by~(\ref{11}).
	Suppose
\begin{itemize}
  \setlength{\itemsep}{0cm}%
  \setlength{\parskip}{0cm}%

      \item[$\star$]  $\w\phi$ is the distribution associated with $L$;

 	\item[$\star$] $\phi \in {\cal L}_p $ and
	there exists $B_{\phi}>0$ such that
$
\sum\limits_{k\in\zd}  |\h\phi(\xi+k)|^q<B_{\phi}\quad\forall \xi\in\rd;
$
there exist $n\in\n$ and $\delta\in(0, 1/2)$ such that
	$\h\phi\h{\w\phi}$ is  boundedly differentiable up to order $n$ on
	$\{|\xi|<\delta\}$,
	$\h\phi$ is boundedly	differentiable up to order $n$  on $\{|\xi+l|<\delta\}$	for all $l\in\zd\setminus\{\nul\}$;
	the function $\sum\limits_{l\in\zd,\, l\neq\nul}|D^\beta \h \phi (\xi+l)|$
is bounded on	$\{|\xi|<\delta\}$  for $[\beta]=n$;
	$D^{\beta}(1-\h\phi\h{\w\phi})(0) = 0$ for all $\beta\in\Delta_{n}$;
	the Strang-Fix condition of order $n$ holds 	for $\phi$;
	
      \item[$\star$]  $f=\h g^{-}$, where
       $g\in L_q$, $g(\xi)=O(|\xi|^{-N-d-\varepsilon})$
as $|\xi|\to\infty$,
	$\varepsilon>0$.

\end{itemize}

\noindent
	Then
	\be
	\left\|f-m^{-j/2} \sum_{k\in\zd}L f(M^{-j}\cdot)(-k) \phi_{jk}\right\|^q_p\le
	C_1
	 \|M^{*-j}\|^{\gamma q}  {\cal I}_{j,\gamma,q}^{Out}(g)+
	 C_2 \|M^{*-j}\|^{nq}   {\cal I}_{j,n,q}^{In}(g),
	 \label{fTheoQjMain1}
	\ee
where $C_1$ and $C_2$ do not depend on $j$ and $f$.

If, in addition, $n< N+\frac dp + \varepsilon$, then
$$
 \left\|f-m^{-j/2}\sum_{k\in\zd} L f(M^{-j}\cdot)(-k) \phi_{jk}\right\|_p\le C
 \theta^{-jn}
 $$
  for every  positive number $\theta$   which is smaller in module
than any eigenvalue of  $M$ and  some $C$ which does not depend on $j$.

If, in addition, $M$ is an isotropic matrix dilation and $\lambda$
is its eigenvalue then
	 \be
	 \left\|f-m^{-j/2}\sum_{k\in\zd} L f(M^{-j}\cdot)(-k) \phi_{jk}\right\|_p\le
	 \begin{cases}
	 C |\lambda|^{-j(N+\frac dp + \varepsilon)}  &\mbox{if }
	n> N+\frac dp + \varepsilon\\
	  C j^{1/q} |\lambda|^{-jn} &\mbox{if }
	 n= N+\frac dp + \varepsilon \\
	C|\lambda|^{-jn}
	 &\mbox{if }
	 n< N+\frac dp + \varepsilon
	\end{cases},
	\label{fTheoQjL}
	\ee
	where $C$ does not depends on $j$.
\end{theo}	

Similarly, the following statement follows from Theorem~\ref{theoQj1}.

 	\begin{theo}
\label{theoQj11}
	Let $2\le p \le \infty$, $1/p+1/q=1$, a differential operator $L$
be defined by~(\ref{11}).
	Suppose
\begin{itemize}
  \setlength{\itemsep}{0cm}%
  \setlength{\parskip}{0cm}%

      \item[$\star$]  $\w\phi$ is the distribution associated with $L$;

\item[$\star$] $\phi \in {\cal L}_p $, the function
 $\sum_{k\in\zd}  |\h\phi(\xi+k)|^q$	is bounded, and
     there exists $\delta\in(0,1/2)$ such that
 $\overline{\h\phi}(\xi)\h{\w\phi}(\xi)=1$
 a. e. on $\{|\xi|<\delta\}$, $\h\phi(\xi)=0$ a.e.
 on $\{|l-\xi|<\delta\}$ for all
 $l\in\z\setminus\{0\}$;

      \item[$\star$]  $f=\h g^{-}$, where
       $g\in L_q$, $g(\xi)=O(|\xi|^{-N-d-\varepsilon})$
as $|\xi|\to\infty$,  $\varepsilon>0$.

\end{itemize}

\noindent
	Then
 \be
 \left\|f-m^{-j/2}\sum_{k\in\zd} L f(M^{-j}\cdot)(-k) \phi_{jk}\right\|_p\le C
 \theta^{-j(N+\frac{d}{p}+\varepsilon)}
 \label{12}
 \ee
 for every  positive number $\theta$   which is smaller in module
than any eigenvalue of  $M$ and  some $C$ which does not depend on $j$.
\end{theo}	

\section{Falsified sampling expansions}
\label{FalSamExp}

 If exact sampled values of a signal $f$ are  known
 then sampling expansions are very useful for applications.
 Theorems~\ref{theoQjL} and \ref{theoQj11} provide error estimates
 for this case.
  We now discuss what happens
 if exact sampled values at the points $M^{-j}k$
  are replaced by the following average values 
\be   
	\frac  1 {V_h} \int\limits_{B_h} f(M^{-j}k+M^{-j}t)d t=
	\frac  {m^j} {V_h} \int\limits_{M^{-j}B_h} f(M^{-j}k+t)d t=: 
{\rm \bf Av}_{h}(f, M^{-j}k),
	\label{fFAverValue}
	\ee
where $V_h$ is the volume of the ball $B_h,$  $h>0.$ In the case $d=1$, $M=2$,
it easily follows from the Tailor formula (see Introduction) that
$$
\frac1{2^{-j}h}\int\limits_{2^{-j}k}^{2^{-j}(k+h)}f(t)\,dt\approx
Lf(2^{-j}\cdot)(k),
$$
where $L$ is the differential operator~(\ref{11})
with $\alpha_l=\frac1{(l+1)!}h^l$.

We need the following lemma 	
to prove a similar statement for the averages~(\ref{fFAverValue}).
\begin{lem}
Let $N\in\n,$ function $f$
be continuously differentiable up to order $N$, $A$ be a real-valued $d\times d$ matrix.
	Then for all $t,x\in\rd$
	$$
	\sum\limits_{\beta\in\Delta_{N+1}}
	\frac{D^{\beta} f (A x)} {\beta!}
	(A t)^{\beta}=
	\sum\limits_{\beta\in\Delta_{N+1}}
	\frac{ D^{\beta} [f(Ax)]}{\beta!}
	t^{\beta}.
	$$
\label{lemDiffMatrixA}
\end{lem}

{\bf Proof. }
Firstly, we introduce some additional notations. Let $p\in\z_+,$
$O_p=\{\beta\in\zd_+: [\beta]=p\}$. Assume that the set $O_p$ is
ordered by lexicographic order. Namely,
$(\beta_1,\dots,\beta_d)$ is less than
$(\alpha_1,\dots,\alpha_d)$ in lexicographic order
if $\beta_j=\alpha_j$ for $j=1,\dots,i-1$ and
$\beta_i<\alpha_i$ for some $i$.
Let $S(A,p)$ be a  $(\# O_p)\times (\# O_p)$ matrix
 which is uniquely determined by
	\be
	\frac {(A t)^{\alpha}}{\alpha!} = \sum_{\beta\in O_p}
	[S(A,p)]_{\alpha,\beta} \frac {t^{\beta}}{\beta!},
	\label{fBHanEq}	
	\ee
where $\alpha\in O_p,$ $t\in\rd.$ It can be verified that
	\be
	\alpha! [S(A,p)]_{\alpha,\beta}=
	\beta! [S(A^*,p)]_{\beta,\alpha}.
	\label{fBHanEq2}\ee
The above notations and the latter fact was borrowed~\cite{HanSymSmooth}.

Fix $\beta \in\zd_+.$ Let ${\cal E}$ be the set of ordered samples with replacement
of size $[\beta]$ from the set $\{e_1,\dots,e_d\}$, where
$e_k$ is the $k$-th ort in $\rd$. An element $e\in {\cal E} $ is a set  $\{e_{i_1},\dots,e_{i_{[\beta]}}\}$,
where $i_l\in\{1,\dots,d\},$ $l=1,\dots,[\beta],$
$\#{\cal E} = d^{[\beta]}.$ For $e\in {\cal E}$ denote by $(e)_l:=e_{i_l}$,  $l=1,\dots,[\beta].$
Let $T$ be a function defined on ${\cal E}$ by $T(e):=\sum_{i=1}^{[\beta]} (e)_i.$
Note that $T(e)\in \zd_+$ and $[T(e)]=[\beta].$ Denote by $b$ an element of
${\cal E}$ so that $T(b)=\beta$. Such $b$ is unique up to a permutation.
Using the higher chain rule,  we have
	$$
	D^{\beta} [f(Ax)]= \frac {\partial ^{[\beta]} f(A\cdot)}
	{\partial x^{\beta}}(x)=
	\sum_{e\in{\cal E}} \frac {\partial^{[\beta]} f(y)}
	{\partial y^{T(e)}}\Big|_{y=Ax}
	\prod_{i=1}^{[\beta]} \frac {\partial (Ax)^{(e)_i}}
{\partial x^{(b)_i}},\quad x\in\rd,
	$$
	where $\prod_{i=1}^{[\beta]} \frac {\partial (Ax)^{(e)_i}}{\partial x^{(b)_i}}$	 does not depend on $x$ and $D^{\beta} [f(Ax)]$ does not depend on the choice of $b$.
For different elements $e,h\in{\cal E}$, we may have $T(e)=T(h).$
Thus, we can group terms in the sum with equal values of $T(\cdot)$. Namely,
\be
	D^{\beta} [f(Ax)] =
		\sum_{\alpha \in\zd_+, [\alpha]=[\beta]}
		\frac {\partial^{[\beta]} f(y)}
	{\partial y^{\alpha}}\Big|_{y=Ax}
	\sum_{e\in{\cal E}, T(e)=\alpha}
	\prod_{i=1}^{[\beta]} \frac {\partial (Ax)^{(e)_i}}
	{\partial x^{(b)_i}}.
\label{16}
\ee	
If $f(x)=e^{2\pi i (t,x)},$  $t\in\rd,$ then
	$
D^{\beta} [f(Ax)] =	D^{\beta} e^{2\pi i (A^*t,x)}=(A^*t)^{\beta} e^{2\pi i (t,Ax)}.
	$
On the other hand, by~(\ref{16}),
	$$
	D^{\beta} [f(Ax)] =
			 e^{2\pi i (t,Ax)}
			 \sum_{\alpha \in\zd_+, [\alpha]=[\beta]}
t^{\alpha}
	\sum_{e\in{\cal E}, T(e)=\alpha}
	\prod\limits_{i=1}^{[\beta]}
	\frac {\partial (Ax)^{(e)_i}}{\partial x^{(b)_i}}.
	$$
	Thus, due to~(\ref{fBHanEq}) with matrix $A$ replaced by $A^*$,
	 we obtain
	\be\sum_{e\in{\cal E}, T(e)=\alpha}
	\prod\limits_{i=1}^{[\beta]}
	\frac {\partial (Ax)^{(e)_i}}{\partial x^{(b)_i}}=[S(A^*,p)]_{\beta,\alpha} \frac {\beta!}{\alpha!}.
	\label{fLemBHDiff}	
	\ee
	
Let now $f$ be an arbitrary  function continuously differentiable  up
	to order $N,$ $0\le p \le N$, $p\in\z_+.$
	It follows from~(\ref{fBHanEq2}) 	and~(\ref{fLemBHDiff}) that
	
	$$
	\sum_{\alpha\in O_p}
	\frac {D^{\alpha} f(A x)}{\alpha!} (A t)^{\alpha}=
		\sum_{\alpha\in O_p}
	D^{\alpha} f(A x) \sum_{\beta\in O_p}
[S(A,p)]_{\alpha,\beta} \frac {t^{\beta}}{\beta!}=
		\sum_{\beta\in O_p}  t^{\beta}
		\sum_{\alpha\in O_p}
	D^{\alpha} f(A x)
 \frac {[S(A^*,p)]_{\beta,\alpha}}{\alpha!}=
	$$
	$$
	\sum_{\beta \in O_p}  \frac {t^{\beta}}{\beta!}
		\sum_{\alpha\in O_p}
	D^{\alpha} f(A x)
 \sum_{e\in{\cal E}, T(e)=\alpha}
	\prod\limits_{i=1}^{[\beta]} \frac {\partial (Ax)^{(e)_i}}
	{\partial x^{(b)_i}}=
	\sum_{\beta\in O_p}  \frac {t^{\beta}}{\beta!}
	 D^{\beta} [f(Ax)].
	$$
It remains to sum the latter expression over $p$ from $0$ to $N.$
$\Diamond$

Let
	\be
	L f(M^{-j}\cdot)(k)=
	\sum_{\beta\in\Delta_{N+1}} a_\beta D^{\beta}
	[f(M^{-j}\cdot)](k) ,
	\quad
	a_\beta =\frac  {1} {\beta! V_h} \int\limits_{B_h} t^{\beta}d t.
	\label{fOperL}
	\ee
By Lemma~\ref{lemDiffMatrixA}, we have
	\ba
\nonumber
	L f(M^{-j}\cdot)(k)=
	 \frac  {1} { V_h} \int\limits_{B_h}
	 \sum_{\beta\in\Delta_{N+1}}	
	\frac {D^{\beta}
	[f(M^{-j}\cdot)](k)}{\beta!}
	t^{\beta}d t =
	 	\\
\frac  {1} { V_h} \int\limits_{B_h}
	 \sum_{\beta\in\Delta_{N+1}}	
	\frac{D^{\beta} f (M^{-j} k)} {\beta!}
	(M^{-j} t)^{\beta} d t =
\frac {m^j}{V_h}
	\int\limits_{M^{-j}B_h}
	 \sum_{\beta\in\Delta_{N+1}}	
	\frac{D^{\beta} f (M^{-j} k)} {\beta!} 		t^{\beta}d t,
	\label{15}
\ea
and, due to the Tailor formula,
$$
\frac  {m^{j}} {V_h} \int\limits_{M^{-j}B_h} f(M^{-j}k+t)d t
	\approx L f(M^{-j})(k).
$$

We now are interested in error analysis for {\em falsified sampling expansions}
$$
m^{-j/2}\sum_{k\in\zd}{\rm \bf Av}_{h}(f, M^{-j}k)\,\phi_{jk}
=\sum_{k\in\zd}
\frac  {m^{j/2}} {V_h} \int\limits_{M^{-j}B_h} f(M^{-j}k+t)d t\,\phi_{jk}.
$$
To investigate the convergence and  approximation order
of falsified sampling expansions we can use
Theorem~\ref{theoQjL} or Theorem~\ref{theoQj11},   and estimate the sum
$
m^{-j/2}\sum_{k\in\zd}\varepsilon_j(-k)\,\phi_{jk},
$
where
	\ba
	\varepsilon_j(k):=
	\frac  {m^j} {V_h} \int\limits_{M^{-j}B_h} f(M^{-j}k+t)d t-
	L f(M^{-j}\cdot)(k).
		\label{fEpsK}
	\ea

\begin{theo}
Let $d< p \le \infty$,  $h>0$, $N\in\n$, $\phi \in {\cal L}_p$.
Suppose  $f\in W_p^{N+1}$, operator $L$ is defined by~(\ref{fOperL}),
$\varepsilon_j(k)$ is defined by~(\ref{fEpsK}). Then
\be
	\left\| m^{-j/2}
	\sum\limits_{k\in\zd}
	\varepsilon_j(-k) \phi_{jk}\right\|_{p} \le
	C \theta^{-j(N+1)}.
\label{13}
\ee
for every  positive number $\theta$   which is smaller in module
than any eigenvalue of  $M$ and  some  $C$ which  does not depend on $j$.
\label{t1}
\end{theo}

{\bf Proof.} Let fix $j\in\n.$	
Due to Proposition~\ref{propLp},
	\be
	\left\| m^{-j/2}	\sum\limits_{k\in\zd}
	\varepsilon_j(-k) \phi_{jk}\right\|_{p} =
	m^{-j/p} \left\| \sum\limits_{k\in\zd}
	\varepsilon_j(-k)\phi_{0k}\right\|_p\le
m^{-j/p}\|\phi\|_{{\cal L}_p}\lll\sum\limits_{k\in \zd}|\varepsilon_j(k)|^p\rrr^{1/p}.
	\label{14}
\ee
By the Taylor formula with 	integral remainder, we have
	$$
	f(M^{-j}k+t)= \sum\limits_{\beta\in\Delta_{N+1}}
	\frac {D^{\beta} f(M^{-j}k)}{\beta!} t^{\beta} +
	\sum\limits_{\beta\in\zd_+, [\beta]=N+1}
	\frac {N+1}{\beta!} t^{\beta}
	\int\limits_0^1 (1-\tau)^N D^{\beta} f(M^{-j}k+t\tau)d\tau.
	$$
It follows from~(\ref{15})  that
$$
\varepsilon_j(k)=\frac  {m^j} {V_h} \int\limits_{M^{-j}B_h}\,dt
\sum\limits_{\beta\in\zd_+, [\beta]=N+1}
	\frac {N+1}{\beta!} t^{\beta}
	\int\limits_0^1 (1-\tau)^N D^{\beta} f(M^{-j}k+t\tau)\,d\tau.
$$
	Hence, taking into account that $|t^{\beta}|\le |t|^{[\beta]}$, we get
$$
	|\varepsilon_j(k)|\le
	\sum\limits_{\beta\in\zd_+, [\beta]=N+1}
	\frac {N+1}{\beta!}
	\frac  {m^j} {V_h} \int\limits_{M^{-j}B_h}|t|^{[\beta]}  dt
	\int\limits_0^1  |D^{\beta} f(M^{-j}k+t\tau)|d\tau \le
	$$
	$$
	\sum\limits_{\beta\in\zd_+, [\beta]=N+1}
	\frac {N+1}{\beta!} \frac  {m^j} {V_h} 	\int\limits_0^1   d\tau
	 \int\limits_{\tau M^{-j}B_h}\frac {|t|^{[\beta]} }{\tau^{[\beta]}}
	| D^{\beta} f(M^{-j}k+t) |\frac {dt}{\tau^d} =
	 $$
	 $$
	 \sum\limits_{\beta\in\zd_+, [\beta]=N+1}
	 \frac {N+1}{\beta! V_h} 	\int\limits_0^1   d\tau
	 \int\limits_{\tau B_h}|M^{-j}t|^{[\beta]}
	| D^{\beta} f(M^{-j}k+M^{-j}t)| \frac {dt}{\tau^{[\beta]+d}}.
	 $$
 The latter integration is taken over the set
 $\{\tau\in[0,1], t\in \tau B_h\}=\{\tau\in[0,1], |t|\le \tau h \}$, or equivalently
$\{|t|\in[0,h], \frac {|t|}{h}\le\tau \le 1\}$. Changing the order of integration,
we obtain
	 $$|\varepsilon_j(k)|\le
	 \sum\limits_{\beta\in\zd_+, [\beta]=N+1}
	 \frac {N+1}{\beta! V_h} \int\limits_{B_h} |M^{-j}t|^{[\beta]}
	 |D^{\beta} f(M^{-j}k+M^{-j}t)| dt \int\limits_{\frac {|t|}{h}}^1
	\frac {d\tau}  {\tau^{N+d+1}} \le
	$$
	$$
	\sum\limits_{\beta\in\zd_+, [\beta]=N+1}
	 \frac {2}{\beta! V_h}
	\int\limits_{B_h}
	 |M^{-j}t|^{N+1}
	| D^{\beta} f(M^{-j}k+M^{-j}t) |
	  \left(\frac {h}{|t|}\right)^{N+d} dt \le
	  $$
	  $$
	 	\sum\limits_{\beta\in\zd_+, [\beta]=N+1}
	 \frac {2 h^{N+d} \|M^{-j}\|^{N+1}}{\beta! V_h}
	\int\limits_{B_h}		\frac{| D^{\beta} f(M^{-j}k+M^{-j}t) |}{|t|^{d-1}}
	   dt =
$$
$$
\sum\limits_{\beta\in\zd_+, [\beta]=N+1}
	 \frac {2 h^{N+d}m^j \|M^{-j}\|^{N+1}}{\beta! V_h}
	\int\limits_{M^{-j}B_h}\frac{| D^{\beta} f(M^{-j}k+t) |}{|M^j t|^{d-1}}dt.
	  $$
Using H\"olders's inequality, we have	
$$
\int\limits_{M^{-j}B_h}\frac{| D^{\beta} f(M^{-j}k+t) |}{|M^j t|^{d-1}}dt\le
\lll\int\limits_{M^{-j}B_h}{| D^{\beta} f(M^{-j}k+t) |^p}dt\rrr^{1/p}
\lll\int\limits_{M^{-j}B_h}\frac{dt}{|M^j t|^{q(d-1)}}\rrr^{1/q},
$$
where $q=\frac p{p-1}$.
Since
$$
\int\limits_{M^{-j}B_h}\frac{dt}{|M^j t|^{q(d-1)}}\le m^{-j}
\int\limits_{B_h}\frac{dt}{|t|^{q(d-1)}}
$$
and $q(d-1)<d$, the latter integral is finite.
Summarizing the above estimates
we obtain
$$
|\varepsilon_j(k)|^p\le C_1 m^{j} \|M^{-j}\|^{p(N+1)}
\int\limits_{M^{-j}k+M^{-j}B_h}{\sum\limits_{\beta\in\zd_+, [\beta]=N+1}| D^{\beta} f(t) |^p}dt
$$
where $C_1$ does not depend on $f$ and $j.$
It follows that
$$
\sum_{k\in\zd}|\varepsilon_j(k)|^p\le C_1 m^{j} \|M^{-j}\|^{p(N+1)}
\|f\|_{W_p^{N+1}}^p.
$$
Combining this with~(\ref{14}) and~(\ref{00}), we get~(\ref{13}).  $\Diamond$

\begin{rem}
If  $M$ is an isotropic matrix
for which $\lambda$  is an eigenvalue,
then in the proof of Theorem~\ref{t1} we can use inequality~(\ref{10}) instead of~(\ref{00}). Hence inequality~(\ref{13}) can be replaced by
$$
	\left\| m^{-j/2}
	\sum\limits_{k\in\zd}
	\varepsilon_j(-k) \phi_{jk}\right\|_{p} \le
	C |\lambda|^{-j(N+1)}.
$$
\end{rem}

Using the above theorem we can state the convergence and  approximation order of 
falsified sampling expansions.

\begin{theo}
\label{theoQjAverage}
Let   $d< p \le \infty$,  $1/p+1/q=1$, $h>0$, $N\in\z_+,$
$M$ be an isotropic matrix dilation and $\lambda$ be its eigenvalue.
Suppose $\w\phi$ is the distribution associated with the differential
operator $L$ given by~(\ref{fOperL}), $\phi$ and $n$ are as in
Theorem~\ref{theoQjL}; $f=\h g^{-}$, where
       $g\in L_q$, $g(\xi)=O(|\xi|^{-N-d-\varepsilon})$
as $|\xi|\to\infty$,
	$\varepsilon>1$.
Then
\be
\Big\|f-m^{-j/2}\sum_{k\in\zd}{\rm \bf Av}_{h}(f, M^{-j}k)\,\phi_{jk}\Big\|_p\le
 \begin{cases}
	 	C |\lambda|^{-j(N+1)}  &\mbox{if }
	
	n> N+1
\\
	C|\lambda|^{-jn}
	 &\mbox{if }
	 n\le N+1
	 	 \\
	\end{cases},
	 \label{fTheoQj3}
	\ee
where $C$ does not depend on $j$.
\end{theo}

\begin{theo}
\label{theoQjAverage2}
Let   $d< p \le \infty$,  $1/p+1/q=1$, $h>0$, $N\in\z_+.$
Suppose $\w\phi$ is the distribution associated with the differential
operator $L$ given by~(\ref{fOperL}), $\phi$ is as in
Theorem~\ref{theoQj11}; $f=\h g^{-}$, where
       $g\in L_q$, $g(\xi)=O(|\xi|^{-N-d-\varepsilon})$
as $|\xi|\to\infty$,
	$\varepsilon>1$.
Then
\be
\Big\|f-m^{-j/2}\sum_{k\in\zd}{\rm \bf Av}_{h}(f, M^{-j}k)\,\phi_{jk}\Big\|_p\le C \theta^{-j(N+1)},
 	 \label{fTheoQj32}
	\ee
for every  positive number $\theta$   which is smaller in module
than any eigenvalue of  $M$ and  some  $C$ which  does not depend on $j$.
\end{theo}

Observing the proof of Theorem~\ref{t1}, one can see that in the
one-dimensional case  an analog of~(\ref{13})
holds  true for    a wider class of functions $f$ and any $p\ge1$.
Indeed, in this case $M$ is a dilation factor, let $M>0.$
Then we have
$$
|\varepsilon(k)|\le
	 \frac {2}{(N+1)! V_h}
	\int\limits_{B_h}
	 |M^{-j}t|^{N+1}
	|  f^{(N+1)}(M^{-j}k+M^{-j}t) |
	  \left(\frac {h}{|t|}\right)^{N+1} dt \le
$$
$$
\frac {2 h^{N+1}M^{-j(N+1)}}{(N+1)! V_h}
	\int\limits_{B_h}	 	|  f^{(N+1)}(M^{-j}k+M^{-j}t) |	   dt =
\frac {2 h^{N+1}M^{-jN}}{(N+1)! V_h}
	\int\limits_{M^{-j}B_h+M^{-j}k}	|  f^{(N+1)}(t) |dt.
$$
It follows that
$$
\left\| m^{-j/2}	\sum\limits_{k\in\z}
	\varepsilon(-k) \phi_{jk}\right\|_{p} \le
	m^{-j/p}\|\phi\|_p
\sum_{k\in\z}|\varepsilon(k)|\le C_1\|f\|_{W_{1}^{N+1}}M^{-j(N+\frac1p)},
$$
where $C_1$ does not depend on $f$ and $j.$ This yields the following statements.

\begin{theo}
Let $d=1$, $p\ge1$,  $h>0$, $N\in\n$, $\phi\in L_p$,
Suppose  $f\in W_1^{N+1}$,
$\varepsilon(k)$ is defined by~(\ref{fEpsK}). Then
\be
	\left\| M^{-j/2}
	\sum\limits_{k\in\z}
	\varepsilon(-k) \phi_{jk}\right\|_{p} \le  C M^{-j(N+\frac 1p)}.
\label{17}
\ee
where  $C$ does not depend on $j$.
\label{t11}
\end{theo}

\begin{theo}
\label{theoQjAverage1}
Let $d=1$,  $2\le p \le \infty$,   $h>0$, $N\in\z_+,$.
Suppose $\w\phi$ is the distribution associated with a differential
operator $L$ given by~(\ref{fOperL}), $\phi$ and $n$ are as in
Theorem~\ref{theoQjL}, $f\in W_1^{N+1}$, $f^{(N+1)}\in\ Lip \,\varepsilon$,
$\varepsilon>0$. Then
\be
\Big\|f-m^{-j/2}\sum_{k\in\zd}{\rm \bf Av}_{h}(f, M^{-j}k)\,\phi_{jk}\Big\|_p\le
 \begin{cases}
	 	C |\lambda|^{-j(N+\frac 1p )}  &\mbox{if }
	
	n> N+\frac 1p
\\
	C|\lambda|^{-jn}
	 &\mbox{if }
	 n\le N+\frac 1p
	 	 \\
	\end{cases},
	 \label{fTheoQj31}
	\ee
where $C$ does not depend on $j$.
\end{theo}

	\section{Examples}
	\label{example}

In this section some examples will be given to illustrate the obtained results.

Firstly we discuss construction of band-limited functions $\phi$.
Theoretically, for every differential operator $L$ one can easily
construct compactly supported $\h\phi\in C^{d}(\rd)$  such that
$\h\phi \h{\w\phi}=1$ on a small neighborhood of zero.
All assumptions of Theorem~\ref{theoQj11} are satisfied.
The  approximation order of the corresponding
expansions depends on  how smooth is the function $f$. However such
expansions are not good from the computational point of view.
We will not able to derive explicit formulas for $\phi$ which is needed
for implementations.

 For arbitrary  differential operator $L$ and arbitrary $n$, we can
 construct $\phi$ satisfying all assumptions of Theorem~\ref{theoQjL} as follows.
 Let $\w\phi$ be the distribution associated with $L$, $T$ be  a trigonometric polynomial  such that $T$ and all its derivatives  up to order $d+1$ vanish on the
 boundary of  $\frac12\td$ and $D^{\beta}(1-T \h{\w\phi})(\nul)=0$,  $[\beta]<n$.
 Define $\h\phi$ as the restriction of $T$ onto $\frac12\td$. The  approximation order of corresponding expansions depends on  how smooth is the function $f$ but cannot be better than $n$. Deriving explicit formulas for $\phi$ is possible
 in this case, but they will be too bulky. Probably such a way is also not appropriate for applications.

We now present examples which can be useful in practice.

I. For sampling expansions we can take $\w\phi=\delta$ and
$\phi(x)=\prod_{k=1}^d\left(\frac {\sin  {\pi x_k} } { {\pi x_k} }\right)^2$.
The corresponding  expansion of $f$  interpolates $f$ at the points $M^{-j}k$, $k\in\z^d$.
Since $\h\phi(\xi)=\prod_{k=1}^d(1-|\xi_k|)\chi_{[-1,1]^d}(\xi)$,
the function $\sum_{k\in\zd}  |\h\phi(\xi+k)|^q$ is bounded and $\phi\in  {\cal L}_p$. The Strang-Fix condition for $\phi$ of order 1 is valid.
So, all assumptions of Theorem~\ref{theoQjL} are satisfied, but
 the  approximation order  cannot be better than $n=1$ by this theorem.
  $\Diamond$

II. Now we illustrate Theorem~\ref{theoQj11} by improving the previous example.
Let $\psi(x)=\prod_{k=1}^d\left(\frac {\sin  {\pi x_k} } { {\pi x_k} }\right)^2.$
Define $\h\phi(\xi)=2\h\psi(2\xi)-\h\psi(4\xi).$ Thus,
$\h\phi$ is continuous, has its support inside $\td$ and
equal to $1$ for all $\xi\in [1/4,1/4]^d,$ the function $\sum_{k\in\zd}  |\h\phi(\xi+k)|^q$ is bounded.  Clearly, $\phi(x)=\frac 1 {2^{d-1}}\psi (\frac x2) - \frac 1{4^d}\psi (\frac x4)$ and therefore  $\phi\in  {\cal L}_p$.
Taking again $\w\phi=\delta$ we obtain a sampling expansion with approximation order depending on how smooth is $f$, according to~(\ref{12}).

Note that the same function $\w\phi$ is associated with the differential
operator $Lf=f+\sum_{[\beta]=1}a_\beta D^\beta f$, where $a_\beta =0$.
Hence the functions $\phi, \w\phi$ satisfy all conditions of
Theorem~\ref{theoQjAverage} with $n=2$, $N=1$ and arbitrary $h>0$,
and, according to~(\ref{fTheoQj32}), or~(\ref{fTheoQj3}) in the case of 
 isotropic $M$, the approximation
order of the corresponding  falsified sampling expansions  is $2$
for smooth enough functions $f$.   $\Diamond$

We now are interested in compactly supported functions $\phi$
which have an advantage  for applications
because the corresponding sampling and differential expansions are finite.
	
III. Let
	$$
\h\phi(\xi) =\prod_{k=1}^d\left(\frac {\sin  {\pi x_k} } { {\pi x_k} }\right)^2.
$$
Since $\phi(x)=\prod_{k=1}^d(1-|\xi_k|)\chi_{[-1,1]^d}(\xi),$
$\phi$ is compactly supported and in ${\cal L}_p$.
Also, the function $\sum_{k\in\zd}  |\h\phi(\xi+k)|^q$  is bounded,
$\h\phi$ is continuously differentiable up to any order, the function
$\sum\limits_{l\in\zd, l\neq 0}  |D^{\beta}\h\phi(\xi+l)|$  is bounded
near  the origin for $[\beta]=2.$
Also, the Strang-Fix condition of order~$2$ holds for $\phi$.
The values of $\h\phi$ and its derivatives at the origin are
	$$
	\h\phi(\nul)=1,\quad D^{\beta}\h\phi(\nul)=0,\quad
	[\beta]=1.
	$$
 So, if $\w\phi=\delta,$  then all assumptions
 of Theorem~\ref{theoQjL} are satisfied.
The corresponding sampling expansion of a signal $f$ interpolates
$f$ at the points $M^{-j}k$, $k\in\z^2$, the approximation order
depends on how smooth is  $f$, but, according to~(\ref{fTheoQjL}),
 it cannot be better than $2.$
Again $\w\phi$ is associated with the differential
operator $Lf=f+\sum_{[\beta]=1}a_\beta D^\beta f$, where $a_\beta =0$.
Hence the functions $\phi, \w\phi$ satisfy all conditions of
Theorem~\ref{theoQjAverage} with $n=2$, $N=1$ and arbitrary $h>0$,
and, according to~(\ref{fTheoQj3}), the approximation
order of the corresponding  falsified sampling expansions  is $2$
for smooth enough functions $f$.	$\Diamond$

IV. Let $d=2$,
	$$\h\phi(\xi_1,\xi_2) =
	\frac {1} {(\pi^2\xi_1 \xi_2)^3}
	\left(\sin^3 \pi \xi_1
	\sin^3 \pi \xi_2 +
	b_1\sin^3 \pi \xi_1
	\sin^4 \pi \xi_2+
	b_2\sin^4 \pi \xi_1
	\sin^3 \pi \xi_2\right)$$
Again $\phi$ is compactly supported and in ${\cal L}_p$.
Also, the function $\sum_{k\in\zd}  |\h\phi(\xi+k)|^q$  is bounded,
$\h\phi$ is continuously differentiable up to any order.
Since the trigonometric polynomial in the numerator of $\h\phi$ is bounded, the function
$\sum\limits_{l\in\zd, l\neq 0}  |D^{\beta}\h\phi(\xi+l)|$  is bounded
near  the origin for $[\beta]=3.$
Also, the Strang-Fix condition of order $3$ holds for $\phi$.
The values of $\h\phi$ and its derivatives at the origin are
	$$
	\h\phi(0,0)=1,\quad
	 D^{(1,0)}\h\phi(0,0)=D^{(0,1)}\h\phi(0,0)=0.
	$$
	$$
	 D^{(2,0)}\h\phi(0,0)=\pi^2(2b_1-1),\quad  D^{(0,2)}\h\phi(0,0)=\pi^2(2b_2-1),\quad
	 D^{(1,1)}\h\phi(0,0)=0.
	$$
Now, we choose an appropriate differential operator $L$ in the form $Lf=f+a_{(2,0)} D^{(2,0)}f+a_{(0,2)} D^{(0,2)}f$, or equivalently,
the associated distribution $\w\phi= \delta + \overline{a_{(2,0)}}D^{(2,0)}\delta+
\overline{a_{(0,2)}} D^{(0,2)}\delta$.
Since
	$\h{\w\phi} (\xi)= 1 - 4\pi^2  \overline{a_{(2,0)}}\xi_1
	-4\pi^2  \overline{a_{(0,2)}}\xi_2,$
we have
	$$
\h{\w\phi} (0,0)=1,\quad D^{(2,0)}\h{\w\phi} (0,0) = -4\pi^2 \overline{a_{(2,0)}},	
	\quad D^{(0,2)}\h{\w\phi} (0,0) =  -4\pi^2  \overline{a_{(0,2)}},
$$
	$$
	 D^{(1,0)}\h{\w\phi}(0,0)=D^{(0,1)}\h{\w\phi}(0,0)= D^{(1,1)}\h{\w\phi}(0,0)=0.
	$$
To satisfy condition $D^{\beta}(1-\h\phi \h{\w\phi})(0,0)=0$ for $\beta\in\Delta_3$ we have to provide
	$$
	b_1 =\frac 12 ( 1- 4  \overline{a_{(2,0)}}), \quad
	b_2 = \frac 12 (1 - 4 \overline{a_{(0,2)}}).
	$$
Finally, all conditions of Theorem~\ref{theoQjL} are satisfied.
 The  approximation order depends on how smooth is
 $f$, but, according to~(\ref{fTheoQjL}), it cannot be better than $n=3.$

 We now show that the coefficients $b_1, b_2$ can be chosen
such that all conditions of Theorem~\ref{theoQjAverage} are satisfied
with $n=3$, $N=2$ and arbitrary $h>0$.
In this case the differential operator $L$ is given by~(\ref{fOperL}). The
 coefficients $a_{\beta},$ $\beta\in\Delta_3$ are as follows
	$$ a_{0,0}=1, \quad a_{1,0} = a_{0,1}= a_{1,1}= 0, \quad
	a_{2,0}= a_{0,2} = \frac 18 h^2.$$
	Thus, we set
	  $b_1 = b_2 = \frac 12 (1-\frac 12 h^2).$
According to~(\ref{fTheoQj3}), the approximation
order of the corresponding  falsified sampling expansions  is $3$
for smooth enough functions $f$. $\Diamond$	

V. Let $d=1,$
	$$\h\phi(\xi)=
	\frac {\sin^4 \pi\xi + b_1 \sin^5 \pi \xi +
	b_2 \sin^6 \pi\xi + b_3 \sin^7 \pi\xi}
	{(\pi \xi)^4}.$$
Since $\phi$ is bounded and compactly supported,
it is in ${\cal L}_p$,
and $\h\phi$ is  continuously differentiable up to any order.
Also, $\sum_{k\in\zd}  |\h\phi(\xi+k)|^q$  is bounded.
Since the trigonometric polynomial in the numerator of $\h\phi$ is bounded, the
function $\sum\limits_{l\in\zd, l\neq 0}  |D^{\beta}\h\phi(\xi+l)|$  is bounded
near  the origin for $\beta=4.$
Also the Strang-Fix condition of order $4$ holds for $\phi$.
The values of $\h\phi$ and its derivatives at the origin are
	$$
	\h\phi(0)=1,\quad \h\phi'(0)=b_1 \pi,\quad
	\h\phi''(0)=\frac 23\pi^2 (3 b_2 - 2),\quad
	\h\phi'''(0)=\pi^3 (6 b_3 - 5 b_1)
	$$
Now, we choose the appropriate differential operator $L$ in the form $Lf=f+a_1 f'+a_2 f'' + a_3 f'''$, or equivalently,
the associated distribution $\w\phi= \delta - \overline{a_1}\delta'+
\overline{a_2} \delta'' - \overline{a_3} \delta'''$.
Since
	$\h{\w\phi} (\xi)= 1 - 2\pi i \overline{a_1}\xi
	-4\pi^2 \overline{a_2} \xi^2 + 8 \pi^3 i
	\overline{a_3} \xi^3 ,$
we have
	$$\h{\w\phi} (0)=1,\quad \h{\w\phi}' (0) = -2\pi i  \overline{a_1},	
	\quad \h{\w\phi}'' (0) = -8\pi^2 \overline{a_2},\quad
	\h{\w\phi}''' (0) = 48 \pi^3 i 	\overline{a_3}$$
To satisfy condition $D^{\beta}(1-\h\phi \h{\w\phi})(0)=0$ for $\beta=0,1,2,3$
we have to provide
\ba
(1-\h\phi\h{\w\phi})'(0)=\pi (b_1-2 i \overline{a_1})=0,
\nonumber
\\
(1-\h\phi\h{\w\phi})''(0)=\frac 23 \pi^2 (-2+3 b_2  -12 \overline{a_2} - 6 i \overline{a_1} b_1)=0,
\nonumber
\\
(1-\h\phi\h{\w\phi})'''(0)=-\pi^3 (5 b_1 - 6 b_3 + 4 i (3b_2 - 2 )\overline{a_1} + 24 b_1 \overline{a_2} - 48 i \overline{a_3})=0.
\label{fExamCoef}
\ea
Thus, the coefficients of the function $\h\phi$ can be easily  successively found  using the coefficients of the differential operator $L$.
Finally, all conditions of Theorem~\ref{theoQjL} are satisfied.
 The  approximation order depends on how smooth is
 $f$, but, according to~(\ref{fTheoQjL}) with $|\lambda|=|M|$, it cannot be better than $n=4.$

We now show that the coefficients $b_1, b_2, b_3$ can be chosen
such that all conditions of Theorem~\ref{theoQjAverage} are satisfied with
$n=4$, $N=3$ and arbitrary $h>0$.
In this case the differential operator $Lf=a_0f+a_1 f'+a_2 f'' + a_3 f'''$ is given by~(\ref{fOperL}) and its
 coefficients are defined as
	$$ a_0=1, a_1 = 0, a_2= \frac {h^2} 6, a_3=0.$$
	Using~(\ref{fExamCoef}),
	 we set $b_1 = 0, b_2 =
	-\frac 23 h^2 - \frac 23, b_3 = 0.$
According to~(\ref{fTheoQj3}), the approximation
order of the corresponding  falsified sampling expansions  is $4$
for smooth enough functions $f$.

Note that all conditions of Theorem~\ref{theoQjAverage1} are also satisfied,
which provides approximation order for a wider class of functions $f$.
Namely,  according to~(\ref{fTheoQj31}),
the approximation order is $3+\frac1p$, whenever $f\in W_1^4$, $f^{(IV)}\in Lip\,\varepsilon$, $\varepsilon>0$.  $\Diamond$

\end{document}